\documentclass[reqno,11pt]{amsart}

\usepackage[top=3cm,bottom=2.5cm,left=2.5cm,right=2.5cm]{geometry}

\usepackage{amsmath,amssymb}
\usepackage{graphicx}
\usepackage{enumerate}
\usepackage{color}
\usepackage{hyperref}
\usepackage[numbers,sort&compress]{natbib}

\allowdisplaybreaks

\newtheorem{theorem}{Theorem}[section]
\newtheorem{lemma}[theorem]{Lemma}
\newtheorem{proposition}[theorem]{Proposition}
\newtheorem{corollary}[theorem]{Corollary}

\theoremstyle{definition}
\newtheorem{definition}[theorem]{Definition}
\newtheorem{assump}{Assumption}

\theoremstyle{remark}
\newtheorem{remark}[theorem]{Remark}

\numberwithin{equation}{section}

\newcommand{\R}{\mathbb{R}}

\title{A Study of the PDE Model for the Anisotropic Consensus-based Optimization}

\author{Keyu Li$^{\,1}$, Hui Huang$^{\,1,*}$}

\thanks{H. Huang is partially supported by the Start-up Grant from Hunan University.}
\thanks{Corresponding author: Hui Huang}

\date{}

\begin{document}

\maketitle

\begin{center}
{\footnotesize
$^1$ School of Mathematics, Hunan University, Changsha 410082, P. R. China\\\
Email: huihuang1@hnu.edu.cn; keyuli@hnu.edu.cn\\
}
\end{center}

\begin{abstract}
We investigate the nonlinear Fokker--Planck equation
associated with anisotropic consensus-based optimization.
Under suitable assumptions on the objective function,
we establish global existence and uniqueness of
density-valued weak solutions for bounded initial
probability densities with finite fourth moment.
The main analytical challenge arises from the
component-wise degeneracy of the diffusion, which
vanishes along moving coordinate hyperplanes
determined by the weighted consensus point.
To overcome this difficulty, we develop spatial
translation estimates based on weighted gradient
bounds and uniform density estimates.
Combined with temporal translation and moment
estimates, these yield strong compactness in
$L^2((0,T)\times\mathbb{R}^d)$ through the
Fr\'echet--Kolmogorov compactness criterion,
enabling passage to the limit in the regularized
equations.
Uniqueness is established through the associated
McKean--Vlasov stochastic differential equation
and Wasserstein stability of the weighted
consensus point.
Furthermore, we prove the propagation of
$H^2$-regularity for sufficiently regular
initial data.
\end{abstract}

\noindent
{\small\textbf{MSC (2020):} 35K55, 35K65, 35K67, 35Q70, 35Q93}

\noindent
{\small\textbf{Keywords:} Consensus-based optimization, Component-wise,
Translation estimates, Regular solutions}

\section{Introduction}

Consensus-based optimization (CBO) is a class of
derivative-free, particle-based methods for the global
minimization of possibly nonconvex and nonsmooth
objective functions. The development of CBO is closely
connected to two research directions: population-based
metaheuristic optimization and individual-based models
for collective dynamics. Classical population-based optimization methods include
evolutionary algorithms, such as Genetic Algorithms (GA)
\cite{fogel1995toward,tang1996genetic},
Particle Swarm Optimization (PSO)
\cite{kennedy1995particle}, and Ant Colony Optimization (ACO)
\cite{dorigo2005ant}. Other widely used metaheuristic
approaches include Simulated Annealing (SA)
\cite{kirkpatrick1983optimization}.
These methods aim to balance exploration of the search
space with exploitation of the information gathered
during the optimization process, and have found
applications in image processing
\cite{karaboga2014comprehensive}, machine learning
\cite{atyabi2013pso}, and path planning
\cite{liang2020applying}.
On the other hand, individual-based models have been
extensively studied to describe collective phenomena,
including swarming, crowd dynamics, synchronization,
and opinion formation.

Motivated by these developments, Pinnau et
al.~\cite{pinnau2017consensus} introduced CBO as a
derivative-free optimization method in which
interacting particles evolve toward a weighted
consensus point while stochastic perturbations
promote exploration of the search space.
The relationship between CBO and stochastic
particle swarm optimization has subsequently been
investigated in~\cite{grassi2021particle},
and global convergence results for related PSO
models have been established
in~\cite{huang2023global}.
Since its introduction, CBO has attracted
considerable attention from both numerical and
analytical perspectives.
Numerical studies have demonstrated the competitive
performance of CBO and its variants on various
benchmark problems and practical applications
\cite{beddrich2026constrained,sun2026consensus}.
Meanwhile, a growing body of theoretical work has
established rigorous results on consensus formation
and global or approximate global convergence
under suitable assumptions
\cite{carrillo2018analytical,ha2020convergence,
fornasier2024consensus,huang2026faithful}.

In recent years, a variety of CBO variants and extensions have been
developed to improve its performance and broaden its range of
applications. Some studies aim to enhance exploration, improve
convergence behavior, or reduce computational costs. For example,
memory effects and, when available, gradient information have been
incorporated into CBO~\cite{riedl2024memory}. Momentum-based CBO
methods have been proposed for stochastic optimal
control~\cite{lyu2025control}, while random particle selection has
been investigated as a way to reduce computational
costs~\cite{borghi2023memory}. Other developments address
optimization problems with additional structures. These include
constrained optimization~\cite{beddrich2026constrained} and
nonconvex--nonconcave min--max optimization~\cite{borghi2026minmax}.
CBO has also been adapted to optimization on the
sphere~\cite{fornasier2021consensus}, robust federated
learning~\cite{garcia2025defending}, and trajectory optimization
in robotics~\cite{sun2026consensus}.

We consider the global optimization problem
\begin{equation*}
    \text{Find}\quad
    x_* \in \underset{x\in\R^d}{\arg\min}\, L(x),
\end{equation*}
where $L:\R^d\to\R$ is a possibly nonconvex
and nonsmooth objective function.
To approximate a global minimizer, CBO employs
a system of interacting particles whose dynamics
combine deterministic attraction toward a weighted
consensus point with stochastic perturbations
that promote exploration of the search space.
The consensus point is determined by the objective
function values of the particles, assigning greater
influence to those with lower objective values.
We next recall the corresponding interacting
particle system \cite{pinnau2017consensus,carrillo2018analytical,
totzeck2018numerical}.

A general CBO system of $N$ interacting
particles $\{X_t^{i,N}\}_{i=1}^N$ in $\R^d$ is governed by
\begin{equation}
\label{eq:CBO-particle}
    dX_t^{i,N}
    =
    -\lambda\bigl(
        X_t^{i,N}-m_\beta^L(\rho_t^N)
    \bigr)\,dt
    +
    \sigma D\bigl(
        X_t^{i,N}-m_\beta^L(\rho_t^N)
    \bigr)\,dB_t^i,
    \qquad i=1,\ldots,N,
\end{equation}
where $\lambda>0$ and $\sigma>0$ are the drift and
diffusion parameters, respectively, and
$D:\R^d\to\R^{d\times d}$ denotes the diffusion
matrix. The processes $\{B_t^i\}_{i=1}^N$ are
independent standard $d$-dimensional Brownian motions.
The initial particles $\{X_0^{i,N}\}_{i=1}^N$ are
assumed to be independent and identically distributed
with common law $\rho_0\in\mathcal P_2(\R^d)$.

The empirical measure of the particle system is
defined by $\rho_t^N
    :=
    \frac1N\sum_{j=1}^N\delta_{X_t^{j,N}},$
and the weighted consensus point is given by
\small{\begin{equation}
\label{consensus point}
    m_\beta^L(\rho_t^N)
    :=
    \frac{
        1
    }{
        \displaystyle\sum_{j=1}^N
        e^{-\beta L(X_t^{j,N})}
    }\displaystyle\sum_{j=1}^N
        X_t^{j,N}e^{-\beta L(X_t^{j,N})},
\end{equation}}
where $\beta>0$ controls the influence of the
objective function on the consensus point. The weighted average point (consensus point) plays a crucial role in the
CBO method. The weight function $\omega_\beta^L(x):=e^{-\beta L(x)}$ is determined
jointly by the objective function values and the parameter $\beta$.
More precisely, the smaller the objective function value of a particle is, the
larger its influence becomes; that is, the particle is assigned a larger weight
in the instantaneous weighted average. When $\beta=0$, all particles have the
same weight, whereas as $\beta\to\infty$, $m_\beta^L(\rho_t^N)$ approaches the global
minimizer among all particles, namely,
\[
    m_\beta^L(\rho_t^N)
    \approx
    \arg\min
    \left\{
        L(X_t^1(t)),\cdots,L(X_t^N(t))
    \right\}.
\]
In fact, the motivation for the weight function $m_\beta^L(\rho_t^N)$ originates
from the celebrated Laplace principle~\cite{dembo2009large,miller2006applied}, see also~\cite[Appendix A.2]{wei2025consensus}, a classical asymptotic method
for integrals.

The choice of the diffusion matrix $D$ determines
different variants of the CBO method. In particular,
the classical isotropic CBO model corresponds to
\begin{equation}
\label{eq:isotropic-D}
    D_{\mathrm{iso}}(z)
    =
    |z|I_d,
    \qquad z\in\R^d,
\end{equation}
where $I_d$ is the $d\times d$ identity matrix.
In contrast, the anisotropic, or component-wise,
CBO model introduced in
\cite{carrillo2021consensus}
employs the diagonal diffusion matrix
\begin{equation}
\label{intro:anisotropic-D}
    D_{\mathrm{aniso}}(z)
    =
    \operatorname{diag}(z_1,\ldots,z_d),
    \qquad z=(z_1,\ldots,z_d)\in\R^d.
\end{equation}
In this case, the stochastic perturbation in each
coordinate depends only on the corresponding
coordinate difference from the consensus point.
Compared with isotropic diffusion, this component-wise
structure offers significant advantages for
high-dimensional optimization, including convergence
conditions that can be independent of the ambient
dimension \cite{carrillo2021consensus,
fornasier2022anisotropic}.

More recently, Bonandin, Riedl, and Veneruso
\cite{bonandin2026exploiting} demonstrated that anisotropic CBO can effectively exploit the additively separable structure of objective
functions, with computational complexity governed by the intrinsic dimension of the individual components rather than the ambient dimension.
However, this advantage depends on the alignment between the coordinate directions and the underlying separability structure. Motivated by this observation, Fiedler and Roith
\cite{fiedler2026consensus} proposed a
reparametrization strategy for objective functions that become separable after a suitable linear
change of variables, further extending the applicability of anisotropic CBO.

An important feature of CBO is that the interacting
particle system admits a natural mean-field
description \cite{gerber2025mean,huang2022mean,huang2025uniform,koss2026mean}. Formally, as $N\to\infty$, the empirical
measure $\rho_t^N$ converges to a probability measure
$\rho_t$, which is the law of the nonlinear
McKean--Vlasov process
\begin{equation}
\label{eq:mean-field-SDE}
    dX_t
    =
    -\lambda\bigl(
        X_t-m_\beta^L(\rho_t)
    \bigr)\,dt
    +
    \sigma D\bigl(
        X_t-m_\beta^L(\rho_t)
    \bigr)\,dB_t,
    \qquad
    \rho_t=\operatorname{Law}(X_t),
\end{equation}
where $B_t$ is a standard $d$-dimensional
Brownian motion and
\small{\begin{equation}
\label{eq:consensus-point}
    m_\beta^L(\rho_t)
    :=
    \frac{
        \displaystyle\int_{\R^d}
        x e^{-\beta L(x)}\,d\rho_t(x)
    }{
        \displaystyle\int_{\R^d}
        e^{-\beta L(x)}\,d\rho_t(x)
    }.
\end{equation}}
By applying It\^o's formula to
\eqref{eq:mean-field-SDE}, the probability density
$\rho=\rho(t,x)$, whenever it exists, formally
satisfies the nonlinear, nonlocal Fokker--Planck
equation
\begin{equation}
\label{eq:general-CBO-PDE}
\begin{aligned}
    \partial_t\rho
    =
    \lambda\nabla\cdot
    \bigl(
        (x-m_\beta^L(\rho_t))\rho
    \bigr)
  +
    \frac{\sigma^2}{2}
    \sum_{k,\ell=1}^d
    \partial_{k\ell}
    \bigl(
        a_{k\ell}(x,\rho_t)\rho
    \bigr),
\end{aligned}
\end{equation}
where
\[
    a(x,\rho_t)
    :=
    D\bigl(x-m_\beta^L(\rho_t)\bigr)
    D\bigl(x-m_\beta^L(\rho_t)\bigr)^\top.
\]

For the isotropic diffusion matrix
\eqref{eq:isotropic-D}, equation
\eqref{eq:general-CBO-PDE} reduces to
\begin{equation}
\label{intro:isotropic}
    \partial_t\rho
    =
    \lambda\nabla\cdot
    \bigl(
        (x-m_\beta^L(\rho_t))\rho
    \bigr)
    +
    \frac{\sigma^2}{2}
    \Delta
    \bigl(
        |x-m_\beta^L(\rho_t)|^2\rho
    \bigr).
\end{equation} In the present work, we focus on the anisotropic
diffusion matrix \eqref{intro:anisotropic-D}.
In this case,
\[
    a(x,\rho_t)
    =
    \operatorname{diag}
    \left(
        (x_1-m_{\beta,1}^L(\rho_t))^2,
        \ldots,
        (x_d-m_{\beta,d}^L(\rho_t))^2
    \right).
\]
Consequently, the corresponding mean-field equation
takes the form
\begin{equation}
\label{component CBO}
    \partial_t\rho
    =
    \lambda\nabla\cdot
    \bigl(
        (x-m_\beta^L(\rho_t))\rho
    \bigr)
    +
    \frac{\sigma^2}{2}
    \sum_{i=1}^d
    \partial_{ii}
    \left(
        (x_i-m_{\beta,i}^L(\rho_t))^2\rho
    \right),
\end{equation}
where $m_{\beta,i}^L(\rho_t)$ denotes the $i$-th
component of the weighted consensus point.

Equation~\eqref{component CBO} is nonlinear and
nonlocal due to the dependence of the weighted
consensus point on the probability distribution
$\rho_t$. Moreover, the diffusion is degenerate:
the diffusion coefficient in the $i$-th coordinate
vanishes along the hyperplane
\[
\left\{
x\in\R^d:
x_i=m_{\beta,i}^L(\rho_t)
\right\}.
\]
This degeneracy differs significantly from that
of the isotropic model~\eqref{intro:isotropic},
whose diffusion vanishes only at the consensus
point. The presence of these coordinate-wise
degeneracies introduces additional analytical
difficulties, particularly in establishing the
spatial compactness of approximate solutions. The main objective of this paper is to establish
the global existence and uniqueness of
density-valued weak solutions to the anisotropic
CBO equation~\eqref{component CBO}, together with
higher regularity properties under suitable
assumptions on the initial data.
A central challenge is to develop a compactness
argument that accommodates the component-wise
degeneracy of the diffusion.

The well-posedness theory for the isotropic CBO equation provides a useful point of comparison. Recently, the existence of weak solutions to the isotropic CBO
equation \eqref{intro:isotropic} was established in \cite{fornasier2025regularity,huang2026mathematical,Sun2026Regularization}. 
Fornasier and Sun \cite{fornasier2025regularity} investigated
regularity and positivity properties of the isotropic
equation and their implications for global convergence. In
\cite{huang2026mathematical}, we established the existence and
uniqueness of weak solutions for bounded initial probability
densities with finite fourth moment, and proved propagation
of $H^2$-regularity for sufficiently regular initial data.
In particular, the proof in \cite{huang2026mathematical} relies on a regularization
procedure, uniform a priori estimates, and a compactness argument
based on the Aubin--Lions lemma \cite{simon1986compact}.
Let \(\rho_\varepsilon^{\mathrm{iso}}\) denote the regularized isotropic densities with cutoff radius
\(R = \varepsilon^{-1}\). In dimensions \(d \ge 3\), the weighted gradient estimates lead to a local Sobolev bound of
the form
\[
\|\rho_\varepsilon^{\mathrm{iso}}\|_{L^2(0,T; W^{1,\gamma}(B_r))} \le C_{T,r,\gamma}, \qquad 1 < \gamma < \frac{2d}{d+2},
\]
for sufficiently small \(\varepsilon\), where \(B_r := B_r(0) = \{x \in \mathbb{R}^d : |x| < r\}\) and the constant is independent of
\(\varepsilon\). Together with suitable estimates for the time derivative, such spatial regularity can be used in
compactness arguments.
More recently, Sun \cite{Sun2026Regularization} established
global well-posedness for initial probability measures with
finite fourth moment under suitable weight assumptions,
and proved smoothing after the consensus point starts
moving, as well as full support at positive times for
non-Dirac initial data.
These results advance the solution theory of the isotropic
equation, while the construction of density-valued weak
solutions for the anisotropic equation requires compactness
estimates adapted to its component-wise degeneracy.

Previous studies have established the well-posedness of
measure-valued solutions associated with anisotropic CBO
through mean-field and stochastic differential equation
approaches~\cite{gerber2025mean}.
In contrast, the present work develops a direct PDE
approach to the global existence and uniqueness of
density-valued weak solutions to the anisotropic
Fokker--Planck equation~\eqref{component CBO}.
While our regularization procedure and basic a priori
estimates build upon the framework developed for the
isotropic CBO equation in~\cite{huang2026mathematical},
the component-wise degeneracy of the diffusion
introduces substantial additional difficulties,
particularly in establishing the spatial compactness
of approximate solutions.
To overcome these difficulties, we develop a new
compactness argument based on spatial translation
estimates adapted to the anisotropic diffusion structure.
To the best of our knowledge, this provides the first
direct PDE-based construction of global density-valued
weak solutions for the anisotropic CBO equation.

To explain this difficulty, let $\rho_\varepsilon$ denote the
regularized anisotropic density with $R=\varepsilon^{-1}$, and set
\[
m_i^\varepsilon(t):=
m_{\beta,i}^L(\rho_\varepsilon(t)).
\]
For almost every $t\in(0,T)$ and $1\leq\gamma<2$, a direct
application of H\"older's inequality leads to
\[
\begin{aligned}
\int_{B_r}|\partial_i\rho_\varepsilon|^\gamma\,dx
&\leq
\left(
\int_{B_r}
|x_i-m_i^\varepsilon(t)|^{-\frac{2\gamma}{2-\gamma}}\,dx
\right)^{\frac{2-\gamma}{2}}
\times
\left(
\int_{B_r}
|x_i-m_i^\varepsilon(t)|^2
|\partial_i\rho_\varepsilon|^2\,dx
\right)^{\frac{\gamma}{2}}.
\end{aligned}
\]
The local weighted gradient estimate established in Section~\ref{The translations estimates}
provides uniform control of the weighted term after integration
in time.
However, whenever the hyperplane
$\{x:x_i=m_i^\varepsilon(t)\}$ intersects the interior of $B_r$,
we have
\[
\int_{B_r}
|x_i-m_i^\varepsilon(t)|^{-\frac{2\gamma}{2-\gamma}}\,dx
=\infty,
\]
because $\frac{2\gamma}{2-\gamma}\geq 2>1.$ Hence, this direct H\"older argument does not provide
a uniform local $W^{1,\gamma}$ estimate.

We therefore obtain spatial compactness by directly estimating
spatial translations; see~\cite{ELIAS202139} for related
compactness arguments.
More precisely, for every fixed $T,r>0$, there exists
$\varepsilon_r\in(0,1)$ such that
\begin{equation*}
\lim_{|y|\to0}
\sup_{0<\varepsilon\leq\varepsilon_r}
\int_0^T\int_{B_r}
|\rho_\varepsilon(t,x+y)-\rho_\varepsilon(t,x)|^2
\,dx\,dt=0.
\end{equation*}
An estimate for the time derivative, combined with spatial
mollification, also yields uniform control of temporal translations;
see~\cite{andreianov2011time} for related time compactness tools.

For any sequence $\varepsilon_n\to0$, the restriction
$\varepsilon_n\leq\varepsilon_r$ holds for all sufficiently large
$n$. The remaining finitely many functions satisfy translation
continuity individually.
The Fr\'echet--Kolmogorov compactness criterion therefore yields,
after passing to a subsequence,
\[
\rho_{\varepsilon_n}\to\rho
\qquad\text{strongly in }L^2((0,T)\times B_r)
\]
for every $r>0$.
In Section~~\ref{existence}, we combine these local translation estimates with
a uniform $L^2$-tail estimate to obtain strong compactness on
the whole space.

The main ingredient of the present compactness argument is thus
a direct spatial translation estimate adapted to the
component-wise diffusion, together with temporal translation
and tail estimates.

We first state the assumptions imposed on the cost function $L$.

\begin{assump}\label{ass}
Throughout the paper, we assume that the objective function
$L\in C(\mathbb{R}^d)$ satisfies the following conditions.
	\begin{enumerate}
		\item $L:~\R^d\to \R$ is bounded from below by $L_m:=\min L$, and there exist constants $C_{L}>0, s\geq 0$ such that
		\begin{equation*}
			|L(x)-L(y)|\leq C_L(1+|x|+|y|)^s|x-y|, \quad \forall x,y\in \R^d;
		\end{equation*}
		%for some constants $L_f>0$ and $s\geq 0$. 
		\item there exist constants $c_\ell,c_u,c_0>0$  and $\ell>0$ such that
		\begin{equation*}
			c_\ell(|x|^\ell-c_0)	\leq L(x)-L_m\leq c_u(|x|^\ell+1),\quad \forall x\in \R^d.
		\end{equation*}
	\end{enumerate}   
\end{assump}

Under Assumption~\ref{ass}, we recall the following estimates for the
consensus point from~\cite[Propositions 3.1 and A.4]{gerber2025mean}.

\begin{proposition}\label{consensus point prop}
Suppose that $L$ satisfies Assumption~\ref{ass} and fix $\beta>0$.
For $K>0$, define
\[
\mathcal P_{2,K}(\mathbb R^d)
:=
\left\{
\mu\in\mathcal P_2(\mathbb R^d):
\int_{\mathbb R^d}|x|^2\,d\mu(x)\leq K
\right\}.
\]
Then there exists a constant
$C_m=C_m(K,\beta,L,d)>0$ such that
\begin{equation}
\label{eq:m-W2}
|m_\beta^L(\mu)-m_\beta^L(\nu)|
\leq C_m W_2(\mu,\nu),
\qquad
\mu,\nu\in\mathcal P_{2,K}(\mathbb R^d).
\end{equation}
Moreover, for every $p\geq1$, there exists a constant
$C_p=C_p(p,\beta,L,d)>0$ such that
\begin{equation}
\label{eq:m-moment}
|m_\beta^L(\mu)|^p
\leq C_p\int_{\mathbb R^d}|x|^p\,d\mu(x),
\qquad
\mu\in\mathcal P_p(\mathbb R^d).
\end{equation}
\end{proposition}

We next formulate the notion of weak solutions
to the nonlinear Fokker--Planck equation
associated with the component-wise CBO model.

\begin{definition}[Definition of weak solutions]\label{Definition of weak solutions}
For any fixed $T>0$, $\rho$ is called a weak solution to the problem
\eqref{component CBO} on $(0,T)$ if $\rho$ satisfies

\begin{enumerate}[(i)]
    \item The function $\rho$ is nonnegative and satisfies
\[
\partial_t\rho\in
L^2(0,T;H^{-2}_{\mathrm{loc}}(\mathbb R^d)),
\qquad
|x|^4\rho\in L^\infty(0,T;L^1(\mathbb R^d)),
\]
and
\[
\rho\in L^\infty(0,T;L^p(\mathbb R^d))
\qquad\text{for every }2\leq p<\infty.
\]
The curve $t\mapsto\rho(t,x)\,dx$ admits a narrowly continuous
representative on $[0,T]$, denoted by $\rho_t$, with
\[
\rho_{t=0}(dx)=\rho_0(x)\,dx.
\]
We identify the density $\rho(t,\cdot)$ with its associated
probability measure whenever no confusion arises.

    \item For any
    $\varphi\in C_c^\infty([0,T)\times\mathbb{R}^d)$,
    the following integral equality holds:
    \begin{align}
-\int_{\mathbb{R}^d}
\rho_0(x)\varphi(0,x)\,dx={}&\int_0^T\int_{\mathbb{R}^d}
\rho\,\partial_t\varphi\,dx\,dt
-\lambda\sum_{i=1}^d
\int_0^T\int_{\mathbb{R}^d}
\bigl(x_i-m_{\beta,i}^L(\rho(t))\bigr)
\rho\,\partial_i\varphi\,dx\,dt
\nonumber\\
&+\frac{\sigma^2}{2}\sum_{i=1}^d
\int_0^T\int_{\mathbb{R}^d}
\bigl(x_i-m_{\beta,i}^L(\rho(t))\bigr)^2
\rho\,\partial_{ii}\varphi\,dx\,dt\nonumber.
\end{align}
\end{enumerate}
\end{definition}

We are now in a position to state our main result,
which establishes the existence and uniqueness
of weak solutions to the component-wise CBO equation.

\begin{theorem}\label{thm:main}
Let $L$ satisfy Assumption~\ref{ass}, and let the non-negative initial data $\rho_0 \in L^1(\mathbb{R}^d)\cap L^\infty(\mathbb{R}^d)\cap\mathcal{P}_4(\mathbb{R}^d).$
Then, for any $T>0$, there exists a  weak solution to the model
\eqref{component CBO} in the sense of Definition~\ref{Definition of weak solutions}.
\end{theorem}
%We shall repeatedly use the uniform bounds
%\begin{equation}\|\nabla\!\cdot h_R\|_{L^\infty}+\sum_{i=1}^d\|\partial_{ii}(h_{R,i}^2)\|_{L^\infty}\le C, \label{eq:hR-derivative-bound}\end{equation}where $C$ is independent of $R\ge1$ and of the density $\rho$.
\begin{theorem}[Uniqueness]
Assume that $
L$ satisfies Assumption~\ref{ass}, and let the non-negative
initial data $\rho_0\in L^1(\mathbb{R}^d)\cap L^\infty(\mathbb{R}^d)
\cap\mathcal{P}_4(\mathbb{R}^d).$
Then the weak solution $\rho$ to \eqref{component CBO} is unique.
\end{theorem}

Moreover, for sufficiently regular initial data, we establish the $H^2$-regularity of solutions, thereby obtaining strong solutions in the local $L^2$ sense.
\begin{lemma}\textup{($H^2$-regularity)}
Let $L$ satisfy Assumption~\ref{ass}, and let $\rho_0$
be a nonnegative probability density such that
$\rho_0\in
L^1(\mathbb R^d)\cap L^\infty(\mathbb R^d)
\cap H^2(\mathbb R^d)\cap\mathcal P_4(\mathbb R^d).$
Then the unique weak solution $\rho$ satisfies
\[
\rho\in L^\infty(0,T;H^2(\mathbb R^d)),
\qquad
\partial_t\rho\in
L^\infty(0,T;L^2_{\mathrm{loc}}(\mathbb R^d)).
\]
\end{lemma}
A key tool in our analysis is the following
Fr\'echet--Kolmogorov compactness theorem
(see \cite[Theorem 4.27, p.~113]{brezis2011functional}),
which we state in the form given in
\cite[Theorem 5]{HANCHEOLSEN2010385} and \cite[Theorem 3.1]{ELIAS202139}.

\begin{theorem}[Fr\'echet--Kolmogorov]
\label{thm:FK}
Let $T>0$, $1\le p<\infty$, and let
$\mathcal F$ be a bounded subset of
$L^p((0,T)\times\mathbb R^d)$.
Assume that the following conditions hold:
\begin{enumerate}
\item[(i)]
For any $\eta>0$ and $r>0$, there exists
$\delta>0$ such that
\[
\begin{aligned}
&\|f(t,x+\xi)-f(t,x)\|_
{L^p((0,T)\times B_r)}+
\|f(t+\tau,x)-f(t,x)\|_
{L^p(I_\tau\times B_r)}
<\eta
\end{aligned}
\]
for all $f\in\mathcal F$ and
$|\xi|+|\tau|<\delta$, where
\[
I_\tau=\{t\in(0,T):t+\tau\in(0,T)\}.
\]

\item[(ii)]
For any $\eta>0$, there exists $r>0$
such that
\[
\|f\|_{L^p((0,T)\times (\mathbb R^d\setminus B_r))}<\eta
\]
for all $f\in\mathcal F$.
\end{enumerate}

Then $\mathcal F$ is relatively compact
in $L^p((0,T)\times\mathbb R^d)$.
\end{theorem}

The remainder of the paper is organized as follows.
In Section~\ref{uniform}, we introduce the regularized problem and
establish uniform a priori estimates, including non-negativity,
conservation of mass, $L^2$ and weighted gradient estimates,
moment bounds, and time derivative estimates.
In Section~\ref{The translations estimates}, we first derive uniform $L^p$ bounds and
then establish spatial and temporal translation estimates
for the approximate solutions.
In Section~\ref{existence}, we obtain a uniform $L^2$-tail estimate and
combine it with the translation estimates to establish
strong compactness on the whole space via the
Fr\'echet--Kolmogorov compactness criterion.
Finally, we pass to the limit in the regularized equation and prove
the existence and uniqueness of global weak solutions to \eqref{component CBO}.
We also establish the propagation of $H^2$ regularity for sufficiently
regular initial data.

\section{A Priori Estimates for the Regularized Problem}
\label{uniform}
%\subsection{Regularized equation}
Throughout the approximation argument, let
$0<\varepsilon\le1$ and $R\ge1$.
Choose $\phi\in C_c^\infty([0,\infty))$ such that
\[
0\le\phi\le1,\qquad
\phi=1\ \text{on }[0,1],\qquad
\phi=0\ \text{on }[2,\infty),
\]
and set
\[
\phi_R(z):=\phi(|z|/R),\qquad z\in\mathbb R.
\]
We consider
\begin{equation}
\begin{aligned}
    \partial_t\rho_{\varepsilon,R}
    =&
    \lambda\sum_{i=1}^d
    \partial_i\!\left(
        \bigl(x_i-m_{\beta,i}^L(\rho_{\varepsilon,R})\bigr)
    \phi_R\bigl(x_i-m_{\beta,i}^L(\rho_{\varepsilon,R})\bigr)
        \rho_{\varepsilon,R}
    \right)\nonumber\\&+
    \frac{\sigma^2}{2}\sum_{i=1}^d
    \partial_{ii}\!\left(
        \bigl(|\bigl(x_i-m_{\beta,i}^L(\rho_{\varepsilon,R})\bigr)
    \phi_R\bigl(x_i-m_{\beta,i}^L(\rho_{\varepsilon,R})\bigr)|^2+\varepsilon^2\bigr)
        \rho_{\varepsilon,R}
    \right),
\end{aligned}
    \label{regularized.}
\end{equation}
For a density $\rho_{\varepsilon,R}$, define
\begin{equation}
    h_{R,i}(t,x)
    :=
    \bigl(x_i-m_{\beta,i}^L(\rho_{\varepsilon,R})\bigr)
    \phi_R\bigl(x_i-m_{\beta,i}^L(\rho_{\varepsilon,R})\bigr),
    \label{def-hR}
\end{equation}
with $h_{R,i}[\rho]$ denoting its $i$th component. We consider
\begin{equation}
\begin{aligned}
    \partial_t\rho_{\varepsilon,R}
    =
    \lambda\sum_{i=1}^d
    \partial_i\!\left(
        h_{R,i}(t,x)
        \rho_{\varepsilon,R}
    \right)+
    \frac{\sigma^2}{2}\sum_{i=1}^d
    \partial_{ii}\!\left(
        \bigl(h_{R,i}(t,x)^2+\varepsilon^2\bigr)
        \rho_{\varepsilon,R}
    \right).
\end{aligned}
    \label{regularized}
\end{equation}

Let $J\in C_c^\infty(\mathbb{R}^d)$ satisfy
\begin{align}\label{J}
J\ge0,\qquad
\int_{\mathbb{R}^d}J(x)\,dx=1,\qquad
\operatorname{supp}J\subset B_1,    
\end{align}
Set \(J_\varepsilon(x)=\varepsilon^{-d}J(x/\varepsilon)\). We prescribe the initial condition
\begin{equation}
    \rho_{\varepsilon,R}(0,x)=
    \rho_{0,\varepsilon}(x):=J_\varepsilon*\rho_0(x),\qquad\|\rho_0\|_{L^1(\R^d)}=1.
    \label{regularized-initial}
\end{equation}
The $h_{R,i}$ satisfies the following estimates for all
$i,j,k\in\{1,\ldots,d\}$; see \cite{huang2026mathematical}:
\begin{align*}
\|h_{R,i}\|_{L^\infty((0,T)\times\R^d)}\le CR,
\qquad
\|\partial_{x_j}h_{R,i}\|_{L^\infty((0,T)\times\R^d)}
\le C\delta_{ij},    
\end{align*}
\[
\|\partial_{x_jx_k}h_{R,i}\|_{L^\infty((0,T)\times\R^d)}
\le \frac{C}{R}\delta_{ij}\delta_{ik}.
\]
Moreover, since
\[
\partial_{ii}(h_{R,i}^2)
=
2(\partial_i h_{R,i})^2
+
2h_{R,i}\partial_{ii}h_{R,i},
\]
we obtain
\begin{align}
\|\nabla\cdot h_R\|_{L^\infty((0,T)\times\R^d)}+
\sum_{i=1}^d\|\partial_{ii}(h_{R,i}^2)\|_{L^\infty((0,T)\times\R^d)}
\le C,\label{h-estimate2}
\end{align}
where $C$ is independent of $R$ and $\varepsilon$.

For each fixed $(\varepsilon,R)$, the diffusion matrix is uniformly elliptic.
Accordingly, a standard fixed-point argument, combined with classical parabolic
theory \cite[Chapter 7]{evans2022partial}, yields a smooth solution to \eqref{regularized} on a nontrivial time
interval. The a priori estimates established below prevent finite-time loss of
the relevant norms and hence allow the solution to be continued up to any
prescribed time $T>0$. 

Firstly, we establish the non-negativity and conservation of mass.

\begin{lemma}[Non-negativity]
\label{lem:nonnegativity}
The solution $\rho_{\varepsilon,R}$ to \eqref{regularized} is non-negative.
\end{lemma}

\begin{proof}
Define $\rho_-:=\min\{0,\rho_{\varepsilon,R}\}$
and, for every $t\in(0,T)$, let
\[
A_t:=\{x\in\mathbb{R}^d:
\rho_{\varepsilon,R}(t,x)<0\}.
\]
Since $\rho_{0,\varepsilon}\geq0$, we have
$\rho_-(0,x)=0.$ Multiplying \eqref{regularized} by $\rho_-$ and integrating
over $\mathbb{R}^d$, we obtain
\begin{align}
\int_{\mathbb{R}^d}\rho_-\partial_t\rho_{\varepsilon,R}\,dx
={}&
\lambda\int_{\mathbb{R}^d}
\rho_-\nabla\cdot(h_R\rho_{\varepsilon,R})\,dx
+
\frac{\sigma^2}{2}
\sum_{i=1}^d
\int_{\mathbb{R}^d}
\rho_-
\partial_{ii}
\left(
(h_{R,i}^2+\varepsilon^2)\rho_{\varepsilon,R}
\right)\,dx.
\label{negative-start}
\end{align}
Since $\rho_-=\rho_{\varepsilon,R}$ on $A_t$ and
$\nabla\rho_-=0$ outside $A_t$, the drift term satisfies
\begin{align}
\lambda\int_{\mathbb{R}^d}
\rho_-\nabla\cdot(h_R\rho_{\varepsilon,R})\,dx
&=
-\lambda\int_{A_t}
\nabla\rho_-\cdot h_R\,\rho_-\,dx
\nonumber\\
&=
-\frac{\lambda}{2}
\int_{A_t}
h_R\cdot\nabla(\rho_-^2)\,dx
=
\frac{\lambda}{2}
\int_{A_t}
\nabla\cdot h_R\,
\rho_-^2\,dx.
\label{negative-drift}
\end{align}
For the diffusion term, for every $i=1,\ldots,d$,
integration by parts gives
\begin{align}
\frac{\sigma^2}{2}\int_{\mathbb{R}^d}
\rho_-
\partial_{ii}
\left(
(h_{R,i}^2+\varepsilon^2)\rho_{\varepsilon,R}
\right)dx
&=
-\frac{\sigma^2}{2}
\int_{A_t}
(h_{R,i}^2+\varepsilon^2)
|\partial_i\rho_-|^2\,dx
+\frac{\sigma^2}{4}
\int_{A_t}
\partial_{ii}(h_{R,i}^2)\rho_-^2\,dx.
\label{negative-diffusion}
\end{align}

Combining \eqref{negative-start}, \eqref{negative-drift}, and \eqref{negative-diffusion}, summing over $i=1,\ldots,d$,
we arrive at
\begin{align}
\frac12\frac{d}{dt}\|\rho_-\|_{L^2(\mathbb{R}^d)}^2
&+
\frac{\sigma^2}{2}
\sum_{i=1}^d
\int_{A_t}
(h_{R,i}^2+\varepsilon^2)
|\partial_i\rho_-|^2\,dx
\nonumber\\
&=
\frac{\lambda}{2}
\int_{A_t}
(\nabla\cdot h_R)\rho_-^2\,dx
+
\frac{\sigma^2}{4}
\sum_{i=1}^d
\int_{A_t}
\partial_{ii}(h_{R,i}^2)\rho_-^2\,dx.
\label{negative-energy}
\end{align}
By \eqref{h-estimate2}, combining with the non-negative diffusion term on the left-hand side of \eqref{negative-energy} yields
\[
\frac{d}{dt}\|\rho_-\|_{L^2(\mathbb{R}^d)}^2 
\leq
C(\lambda+\sigma^2)\|\rho_-\|_{L^2(\mathbb{R}^d)}^2.
\]
Hence, by Gr\"onwall's inequality,
\[
\|\rho_-(t)\|_{L^2(\mathbb{R}^d)}^2
\leq
\exp\!\left\{
C(\lambda+\sigma^2)t
\right\}
\|\rho_-(0)\|_{L^2(\mathbb{R}^d)}^2
=0 \qquad \forall ~t\in(0,T).
\]
Therefore,
\[
\rho_{\varepsilon,R}(t,x)\ge0
\qquad\text{for a.e. }(t,x)\in(0,T)\times\mathbb{R}^d.
\]
This completes the proof.
\end{proof}
Furthermore, the mass is conserved in the following sense
\begin{align}\label{mass-conservation}
\int_{\mathbb{R}^d}
\rho_{\varepsilon,R}(t,x)\,dx
=
\int_{\mathbb{R}^d}
\rho_{0,\varepsilon}(x)\,dx
=\int_{\mathbb{R}^d}
\rho_{0}(x)\,dx=1,
\qquad \text{for a.e. }t\in(0,T).    
\end{align}

We next derive uniform \(L^2\) and weighted gradient estimates for the regularized solutions.

\begin{lemma}[$L^2$-estimate]
\label{lem:L2-weighted-gradient}
Assume that $L$ satisfies Assumption~\ref{ass}. Let the non-negative initial data $\rho_0 \in L^1(\mathbb{R}^d)\cap L^2(\mathbb{R}^d)$ and $\rho_{\varepsilon,R}$ be the solution to the regularized problem \eqref{regularized}, with initial data
$\rho_{0,\varepsilon}$ satisfying $
\|\rho_{0,\varepsilon}\|_{L^2(\mathbb R^d)}
<\infty.$
Then, for any $T>0$,
\begin{equation}
\|\rho_{\varepsilon,R}\|_{L^\infty(0,T;L^2(\mathbb{R}^d))}
\le C,
\label{L2-estimate}
\end{equation}
and
\begin{equation}
\sum_{i=1}^d
\int_0^T\int_{\mathbb{R}^d}
h_{R,i}^2
|\partial_i\rho_{\varepsilon,R}|^2
\,dx\,dt
\le C.
\label{weighted-gradient-estimate}
\end{equation}
%Moreover, the regularization yields
%\begin{equation}
%\varepsilon^2
%\sum_{i=1}^d
%\int_0^T\int_{\mathbb{R}^d}
%|\partial_i\rho_{\varepsilon,R}|^2
%\,dx\,dt
%\le C.
%\label{epsilon-gradient-estimate}
%\end{equation}
where $C$ is a constant independent of $R,\varepsilon$.
\end{lemma}

\begin{proof}
The following computations are justified by the regularity of the approximate solutions and a standard spatial cutoff argument.
We test \eqref{regularized} with $\rho_{\varepsilon,R}$ and integrate
over $\mathbb{R}^d$. This gives
\begin{align}
\int_{\mathbb{R}^d}\rho_{\varepsilon,R}\,\partial_t\rho_{\varepsilon,R}\,dx
={}&
\lambda
\int_{\mathbb{R}^d}
\rho_{\varepsilon,R}\,\nabla\cdot(h_R\rho_{\varepsilon,R})\,dx
+
\frac{\sigma^2}{2}
\sum_{i=1}^d
\int_{\mathbb{R}^d}
\rho_{\varepsilon,R}\,
\partial_{ii}
\left(
(h_{R,i}^2+\varepsilon^2)\rho_{\varepsilon,R}
\right)\,dx.
\label{L2-testing}
\end{align}

Next, we estimate each term on the right-hand side separately. For the drift term, using integration by parts, we obtain
\begin{align}
\lambda
\int_{\mathbb{R}^d}
\rho_{\varepsilon,R}\,\nabla\cdot(h_R\rho_{\varepsilon,R})\,dx
&=
-\lambda
\int_{\mathbb{R}^d}
h_R\rho_{\varepsilon,R}\cdot\nabla\rho_{\varepsilon,R}\,dx
\nonumber\\
&=
-\frac{\lambda}{2}
\int_{\mathbb{R}^d}
h_R\cdot\nabla(\rho_{\varepsilon,R})^2\,dx
=
\frac{\lambda}{2}
\int_{\mathbb{R}^d}
(\nabla\cdot h_R)\rho_{\varepsilon,R}^2\,dx.\label{L2-drift}
\end{align}
For the diffusion term, we obtain
\begin{align}
&\quad\frac{\sigma^2}{2}
\sum_{i=1}^d
\int_{\mathbb{R}^d}
\rho_{\varepsilon,R}\,
\partial_{ii}
\left(
(h_{R,i}^2+\varepsilon^2)\rho_{\varepsilon,R}
\right)\,dx
\nonumber\\
&=
-\frac{\sigma^2}{2}
\sum_{i=1}^d
\int_{\mathbb{R}^d}
(h_{R,i}^2+\varepsilon^2)
|\partial_i\rho_{\varepsilon,R}|^2\,dx
+
\frac{\sigma^2}{4}
\sum_{i=1}^d
\int_{\mathbb{R}^d}
\partial_{ii}(h_{R,i}^2)\rho_{\varepsilon,R}^2\,dx.
\label{L2-diffusion}
\end{align}

Combining \eqref{L2-testing}, \eqref{L2-drift},
and \eqref{L2-diffusion}, we arrive at
\begin{align}
\frac12\frac{d}{dt}\|\rho_{\varepsilon,R}\|_{L^2(\mathbb{R}^d)}^2
&+
\frac{\sigma^2}{2}
\sum_{i=1}^d
\int_{\mathbb{R}^d}
(h_{R,i}^2+\varepsilon^2)
|\partial_i\rho_{\varepsilon,R}|^2\,dx
\nonumber\\
&=
\frac{\lambda}{2}
\int_{\mathbb{R}^d}
(\nabla\cdot h_R)\rho_{\varepsilon,R}^2\,dx
+
\frac{\sigma^2}{4}
\sum_{i=1}^d
\int_{\mathbb{R}^d}
\partial_{ii}(h_{R,i}^2)\rho_{\varepsilon,R}^2\,dx.
\label{L2-energy-identity}
\end{align}
Using the coefficient estimate \eqref{h-estimate2}, we obtain
\begin{equation}
\frac{d}{dt}\|\rho_{\varepsilon,R}\|_{L^2(\mathbb{R}^d)}^2
+
\sigma^2
\sum_{i=1}^d
\int_{\mathbb{R}^d}(h_{R,i}^2+\varepsilon^2)
|\partial_i\rho_{\varepsilon,R}|^2\,dx
\le
C\|\rho_{\varepsilon,R}\|_{L^2(\mathbb{R}^d)}^2.
\label{L2-differential}
\end{equation}

\text{Gr\"onwall}'s inequality therefore gives
\begin{equation*}
\|\rho_{\varepsilon,R}(t)\|_{L^2(\mathbb{R}^d)}^2
\le
e^{Ct}\|\rho_{0,\varepsilon}\|_{L^2(\mathbb{R}^d)}^2,
\qquad \forall ~0\le t\le T.
\end{equation*}
Since $\rho_{0,\varepsilon}=J_\varepsilon*\rho_0$ and
$\|J_\varepsilon\|_{L^1}=1$, Young's inequality gives
$\|\rho_{0,\varepsilon}\|_{L^2}\le\|\rho_0\|_{L^2}$. Hence,
\begin{align}
\sup_{0\le t\le T}\|\rho_{\varepsilon,R}(t)\|_{L^2}^2
\le e^{CT}\|\rho_0\|_{L^2}^2.\label{rhoL2}   
\end{align}
This proves \eqref{L2-estimate}.

Integrating \eqref{L2-differential} over $(0,T)$, we obtain
\begin{align*}
&\|\rho_{\varepsilon,R}(T)\|_{L^2}^2
+\sigma^2\sum_{i=1}^d
\int_0^T\int_{\mathbb R^d}
(h_{R,i}^2+\varepsilon^2)
|\partial_i\rho_{\varepsilon,R}|^2\,dx\,dt\le
\|\rho_{0,\varepsilon}\|_{L^2}^2
+C\int_0^T
\|\rho_{\varepsilon,R}(t)\|_{L^2}^2\,dt.
\end{align*}
Since $\sigma>0$, using \eqref{rhoL2}
and the uniform bound on the initial data,
we conclude that
\begin{equation}
\sum_{i=1}^d
\int_0^T\int_{\mathbb R^d}
(h_{R,i}^2+\varepsilon^2)
|\partial_i\rho_{\varepsilon,R}|^2
\,dx\,dt
\le C,
\label{full-weighted-gradient}
\end{equation}
where $C$ is independent of $R$ and
$\varepsilon$. The proof is complete.
\end{proof}

In addition to the uniform \(L^2\) estimates, we establish moment bounds to control the behavior of the regularized solutions at spatial infinity.

\begin{lemma}[Moment estimates]
\label{lem:moment}
Assume that $L$ satisfies Assumption~\ref{ass}, and let the non-negative initial data $\rho_0 \in L^1(\mathbb{R}^d)\cap L^2(\mathbb{R}^d)\cap\mathcal P_4(\mathbb R^d)$ and $\rho_{\varepsilon,R}$ be the solution to the regularized problem \eqref{regularized}. Then, there exist constants $M_2, M_4> 0$ independent
of $R,\varepsilon$ such that
\begin{align}
&\sup_{0\le t\le T}
\int_{\mathbb R^d}|x|^2\rho_{\varepsilon,R}(t,x)\,dx
\le M_2,
\label{second moment}\\
&\sup_{0\le t\le T}
\int_{\mathbb R^d}|x|^4\rho_{\varepsilon,R}(t,x)\,dx
\le M_4.
\label{fourth moment}
\end{align}  
%where $C$ is a constant independent of $R,\varepsilon$.
\end{lemma}

\begin{proof}
For each fixed $(\varepsilon,R)$, the bounded
coefficients of the regularized SDE imply finiteness of
the second and fourth moments on finite time intervals.
Thus the identities below are justified by smooth
cutoffs and passage to the limit. The estimates derived
from these identities will be uniform in
$\varepsilon$ and $R$.

Since $0\leq \phi_R\leq 1$, we have
\begin{equation}\label{h-est}
|h_{R,i}|
\leq |x_i-m_{\beta,i}^L(\rho_{\varepsilon,R})|,
\qquad
|h_{R,i}|^2
\leq 2|x_i|^2+2|m_{\beta,i}^L(\rho_{\varepsilon,R})|^2.
\end{equation}

We first derive the second moment estimate. Testing the regularized
equation \eqref{regularized} with $|x|^2$ and integrating by parts, we obtain
\begin{align}
\frac{d}{dt}\int_{\mathbb R^d}|x|^2\rho_{\varepsilon,R}\,dx
={}&
-2\lambda\sum_{i=1}^d
\int_{\mathbb R^d}
x_i h_{R,i}\rho_{\varepsilon,R}\,dx
+\sigma^2\sum_{i=1}^d
\int_{\mathbb R^d}
\bigl(h_{R,i}^2+\varepsilon^2\bigr)
\rho_{\varepsilon,R}\,dx.
\label{second moment-eq}
\end{align}
For the first term on the right-hand side of \eqref{second moment-eq}, by \eqref{h-est} and Young's inequality, we have
\[
|x_i h_{R,i}|
\leq |x_i||x_i-m_{\beta,i}^L(\rho_{\varepsilon,R})|
\leq C\bigl(|x_i|^2+|m_{\beta,i}^L(\rho_{\varepsilon,R})|^2\bigr).
\]
Hence,
\begin{align*}
-2\lambda\sum_{i=1}^d
\int_{\mathbb R^d}
x_i h_{R,i}\rho_{\varepsilon,R}\,dx
&\leq
C\sum_{i=1}^d
\int_{\mathbb R^d}
\bigl(|x_i|^2+|m_{\beta,i}^L(\rho_{\varepsilon,R})|^2\bigr)
\rho_{\varepsilon,R}\,dx.
\end{align*}
Since $m_{\beta,i}^L(\rho_{\varepsilon,R})$ is independent of the spatial variable $x$ and
$\rho_{\varepsilon,R}$ has mass conservation \eqref{mass-conservation}, it follows that
\[
\sum_{i=1}^d
\int_{\mathbb R^d}
|m_{\beta,i}^L(\rho_{\varepsilon,R})|^2\rho_{\varepsilon,R}\,dx
=
\sum_{i=1}^d |m_{\beta,i}^L(\rho_{\varepsilon,R})|^2
\int_{\mathbb R^d}\rho_{\varepsilon,R}\,dx
=
|m_{\beta}^L(\rho_{\varepsilon,R})|^2.
\]
Therefore,
\begin{align}\label{L2estimate1}
-2\lambda\sum_{i=1}^d
\int_{\mathbb R^d}
x_i h_{R,i}\rho_{\varepsilon,R}\,dx
\leq
C\left(
|m_{\beta}^L(\rho_{\varepsilon,R})|^2+
\int_{\mathbb R^d}|x|^2\rho_{\varepsilon,R}\,dx
\right).   
\end{align}
On the other hand, since $0<\varepsilon\leq 1$, it follows from
\eqref{h-est} that
\[
h_{R,i}^2+\varepsilon^2
\leq
2\bigl(|x_i|^2+|m_{\beta,i}^L(\rho_{\varepsilon,R})|^2\bigr)+1
\leq
C\bigl(1+|x_i|^2+|m_{\beta,i}^L(\rho_{\varepsilon,R})|^2\bigr).
\]
Consequently, for the second term on the right-hand side of \eqref{second moment-eq}
\begin{align}\label{L2estimate2}
\sigma^2\sum_{i=1}^d
\int_{\mathbb R^d}
\bigl(h_{R,i}^2+\varepsilon^2\bigr)
\rho_{\varepsilon,R}\,dx
&\leq
C\sum_{i=1}^d
\int_{\mathbb R^d}
\bigl(1+|x_i|^2+|m_{\beta,i}^L(\rho_{\varepsilon,R})|^2\bigr)
\rho_{\varepsilon,R}\,dx\nonumber\\
&\leq
C\left(
1+|m_{\beta}^L(\rho_{\varepsilon,R})|^2+
\int_{\mathbb R^d}|x|^2\rho_{\varepsilon,R}\,dx
\right),
\end{align}
where the fixed dimension $d$ has been absorbed into the constant
$C$.

Combining \eqref{L2estimate1}, \eqref{L2estimate2}, and \eqref{second moment-eq}, we obtain
\[
\frac{d}{dt}
\int_{\mathbb R^d}|x|^2\rho_{\varepsilon,R}\,dx
\leq
C\left(
1+|m_{\beta}^L(\rho_{\varepsilon,R})|^2+
\int_{\mathbb R^d}|x|^2\rho_{\varepsilon,R}\,dx
\right).
\]
By Proposition~\ref{consensus point prop}, in particular \eqref{eq:m-moment} with $p=2$, we have $|m_\beta^L(\rho_{\varepsilon,R})|^2
\leq
C\int_{\mathbb R^d}
|x|^2\rho_{\varepsilon,R}\,dx.$
Therefore,
\[
\frac{d}{dt}
\int_{\mathbb R^d}|x|^2\rho_{\varepsilon,R}\,dx
\leq
C\left(
1+
\int_{\mathbb R^d}|x|^2\rho_{\varepsilon,R}\,dx
\right).
\]
Applying Gr\"onwall's inequality once more yields
\[
\sup_{0\leq t\leq T}
\int_{\mathbb R^d}|x|^2
\rho_{\varepsilon,R}(t,x)\,dx
\leq
C
\left(
1+
\int_{\mathbb R^d}|x|^2\rho_{0,\varepsilon}(x)\,dx
\right).
\]
Since $\rho_0\in\mathcal P_2(\mathbb R^d)$,
\[
\sup_{0\leq t\leq T}
\int_{\mathbb R^d}|x|^2
\rho_{\varepsilon,R}(t,x)\,dx
\leq C,
\]
where $C>0$ is independent of $R,\varepsilon$. This proves
\eqref{second moment}.

We next turn to the fourth moment estimate. Testing the regularized
equation with $|x|^4$ and integrating by parts, we obtain
\begin{align}
\frac{d}{dt}
\int_{\mathbb R^d}|x|^4\rho_{\varepsilon,R}\,dx
={}&
-4\lambda\sum_{i=1}^d
\int_{\mathbb R^d}
|x|^2x_i h_{R,i}\rho_{\varepsilon,R}\,dx
\nonumber\\
&
+2\sigma^2\sum_{i=1}^d
\int_{\mathbb R^d}
\bigl(h_{R,i}^2+\varepsilon^2\bigr)
\bigl(|x|^2+2x_i^2\bigr)
\rho_{\varepsilon,R}\,dx.
\label{fourth moment-eq}
\end{align}
For the drift term of \eqref{fourth moment-eq}, by \eqref{h-est},
\begin{equation}\label{fourth-drift}
|x|^2|x_i h_{R,i}|
\leq
|x|^2|x_i|\,|x_i-m_{\beta,i}^L(\rho_{\varepsilon,R})|
\leq
C\bigl(|x|^4+|m_{\beta,i}^L(\rho_{\varepsilon,R})|^4\bigr).
\end{equation}
For the diffusion term, using \eqref{h-est} and
$|x|^2+2x_i^2\leq 3|x|^2$, we have
\begin{align*}
\bigl(h_{R,i}^2+\varepsilon^2\bigr)
\bigl(|x|^2+2x_i^2\bigr)
&\leq
C\bigl(1+|x_i|^2+|m_{\beta,i}^L(\rho_{\varepsilon,R})|^2\bigr)|x|^2\\
&\leq
C\bigl(
|x|^2+|x|^4+|m_{\beta,i}^L(\rho_{\varepsilon,R})|^2|x|^2
\bigr),
\end{align*}
and $|x|^2\leq 1+|x|^4$. Therefore,
\begin{equation}\label{fourth-diffusion}
\bigl(h_{R,i}^2+\varepsilon^2\bigr)
\bigl(|x|^2+2x_i^2\bigr)
\leq
C\bigl(1+|x|^4+|m_{\beta,i}^L(\rho_{\varepsilon,R})|^4\bigr).
\end{equation}

Combining \eqref{fourth moment-eq}--
\eqref{fourth-diffusion} and using mass
conservation, we obtain
\[
\frac{d}{dt}
\int_{\mathbb R^d}|x|^4
\rho_{\varepsilon,R}\,dx
\leq
C\left(
1+
\int_{\mathbb R^d}|x|^4
\rho_{\varepsilon,R}\,dx
+
|m_\beta^L(\rho_{\varepsilon,R})|^4
\right).
\]
By Proposition~\ref{consensus point prop},
in particular \eqref{eq:m-moment} with $p=4$,
we have $|m_\beta^L(\rho_{\varepsilon,R})|^4
\leq
C\int_{\mathbb R^d}
|x|^4\rho_{\varepsilon,R}\,dx.$ Thus,
\[
\frac{d}{dt}
\int_{\mathbb R^d}|x|^4
\rho_{\varepsilon,R}\,dx
\leq
C\left(
1+
\int_{\mathbb R^d}|x|^4
\rho_{\varepsilon,R}\,dx
\right).
\]
Applying Gr\"onwall's inequality once more yields
\[
\sup_{0\leq t\leq T}
\int_{\mathbb R^d}|x|^4
\rho_{\varepsilon,R}(t,x)\,dx
\leq C,
\]
which proves \eqref{fourth moment}. 
\end{proof}

We now turn to the time regularity of the regularized solutions. Using the preceding estimates, we obtain a uniform bound for their time derivatives in a local negative Sobolev space.

\begin{lemma}[Estimate of the time derivative]
\label{time-derivative}
Suppose that $L$ satisfies Assumption~\ref{ass}, and let $\rho_0 \in L^1(\mathbb{R}^d)\cap L^2(\mathbb{R}^d)\cap \mathcal{P}_4(\mathbb{R}^d)$
be non-negative. Let $\rho_{\varepsilon,R}$ be the solution to the
regularized problem \eqref{regularized}. 
Then, for any fixed $T>0$ and bounded open set
$\Omega\Subset\mathbb{R}^d$, there exists a constant
$C_{T,\Omega}>0$, independent of $\varepsilon$ and $R$, such that
\begin{equation}
\|\partial_t\rho_{\varepsilon,R}\|_
{L^\infty(0,T;H^{-2}(\Omega))}
\le C_{T,\Omega}.
\end{equation}
\end{lemma}

\begin{proof}
Fix any bounded open set
$\Omega\Subset\mathbb R^d$.
For any $\psi\in C_c^\infty(\Omega)$,
testing \eqref{regularized} against $\psi$
and integrating by parts, we obtain
\begin{align*}
\left\langle \partial_t\rho_{\varepsilon,R},\psi\right\rangle
={}&
-\lambda\sum_{i=1}^d
\int_{\Omega}h_{R,i}
\rho_{\varepsilon,R}\,\partial_i\psi\,dx
+\frac{\sigma^2}{2}
\sum_{i=1}^d
\int_{\Omega}
\left(h_{R,i}^2+\varepsilon^2\right)
\rho_{\varepsilon,R}\,\partial_{ii}\psi\,dx .
\end{align*}
By (1.11) with $p=2$ and the moment bound \eqref{second moment},
the consensus point is bounded uniformly in $\varepsilon$ and $R$
on $[0,T]$. Hence, since $0\le\phi_R\le1$ and $\Omega$ is bounded,
\[
\|h_{R,i}\|_{L^\infty((0,T)\times\Omega)}
+\|h_{R,i}^2+\varepsilon^2\|_{L^\infty((0,T)\times\Omega)}
\le C_{T,\Omega},
\]
uniformly for $\varepsilon$ and $R$. 
For the first term of the
right side, by the Cauchy-Schwarz inequality,
\[
\begin{aligned}
\left|
\int_{\Omega}
h_{R,i}\rho_{\varepsilon,R}\partial_i\psi\,dx
\right|
&\leq
\|h_{R,i}\|_{L^\infty(\Omega)}
\|\rho_{\varepsilon,R}\|_{L^2(\Omega)}
\|\partial_i\psi\|_{L^2(\Omega)}
\\
&\leq
C_{T,\Omega}
\|\rho_{\varepsilon,R}\|_{L^2(\Omega)}
\|\partial_i\psi\|_{L^2(\Omega)}.
\end{aligned}
\]
Similarly, for the second term,
\[
\left|
\int_{\Omega}
\bigl(h_{R,i}^2+\varepsilon^2\bigr)
\rho_{\varepsilon,R}\partial_{ii}\psi\,dx
\right|
\leq
C_\Omega
\|\rho_{\varepsilon,R}\|_{L^2(\Omega)}
\|\partial_{ii}\psi\|_{L^2(\Omega)}.
\]
Combining these estimates with the uniform
$L^2$ bound \eqref{L2-estimate},
we obtain
\begin{align*}
\left|
\left\langle \partial_t\rho_{\varepsilon,R},\psi\right\rangle
\right|
&\le
C_\Omega
\|\rho_{\varepsilon,R}(t)\|_{L^2(\Omega)}
\left(
\|\nabla\psi\|_{L^2(\Omega)}
+
\|D^2\psi\|_{L^2(\Omega)}
\right)\\
&\le
C_\Omega\|\psi\|_{H^2(\Omega)}.
\end{align*}
The constant $C$ is independent of $\varepsilon,R$.

Since $C_c^\infty(\Omega)$ is dense in $H_0^2(\Omega)$, the above
estimate extends to every $\psi\in H_0^2(\Omega)$. Recalling that
\[
H^{-2}(\Omega)=\bigl(H_0^2(\Omega)\bigr)',
\]
we conclude that
\[
\|\partial_t\rho_{\varepsilon,R}(t)\|_{H^{-2}(\Omega)}
\le C_{T,\Omega}
\]
for almost every $t\in(0,T)$,
uniformly in $\varepsilon$ and $R$.
\end{proof}

\section{The translations estimates}\label{The translations estimates}
In this section, we establish uniform spatial and temporal translation estimates for the approximate solutions to the regularized problem \eqref{regularized}. These estimates are essential for establishing strong compactness via the Fr\'echet--Kolmogorov compactness criterion
(Theorem~\ref{thm:FK}).
%We first derive uniform \(L^p\) bounds to obtain higher integrability of the approximate solutions. We then establish spatial translation estimates, first along each coordinate direction and subsequently for arbitrary spatial translations. Finally, we derive uniform temporal translation estimates.

%In this section, to establish the strong convergence of the approximate solutions, we employ the Fr\'echet--Kolmogorov compactness criterion \ref{thm:FK}. Unlike the isotropic case in \cite{huang2026mathematical}, where compactness can be obtained through local Sobolev regularity estimates by Lions-Aubin lemma \cite{simon1986compact}, the component-wise diffusion considered here requires a different approach. We therefore establish suitable spatial and temporal translation estimates and use them to obtain the relative compactness of the regularized solutions. For completeness, we recall the following compactness criterion 

We now choose
\begin{equation}\label{eq:exist-consensus}
R=\varepsilon^{-1},
\qquad
\rho_\varepsilon:=\rho_{\varepsilon,R},
\end{equation}
and define
\[
h_i^\varepsilon(t,x)
:=
\bigl(x_i-m_{\beta,i}^L(\rho_\varepsilon)\bigr)
\phi_{\varepsilon^{-1}}
\bigl(x_i-m_{\beta,i}^L(\rho_\varepsilon)\bigr).
\]
Then the regularized equation
\eqref{regularized} becomes
\begin{equation}\label{eq:regularized}
\partial_t\rho_\varepsilon
=
\lambda\sum_{i=1}^d
\partial_i\bigl(h_i^\varepsilon(t,x)\rho_\varepsilon\bigr)
+
\frac{\sigma^2}{2}\sum_{i=1}^d
\partial_{ii}
\bigl((|h_i^\varepsilon(t,x)|^2+\varepsilon^2)
\rho_\varepsilon\bigr),
\end{equation}
subject to the initial condition
\begin{equation}
\rho_\varepsilon(0,x)
=
\rho_{0,\varepsilon}(x)
:=
(J_\varepsilon*\rho_0)(x).
\end{equation}

Moreover, since the estimates established in
Lemmas~\ref{lem:nonnegativity}--\ref{time-derivative}
are uniform in both $R$ and $\varepsilon$,
they remain valid for the family
$\{\rho_\varepsilon\}_{0<\varepsilon\le1}$
satisfying \eqref{eq:regularized}.
We summarize these estimates in the following remark.

\begin{remark}
Assume that $L$ satisfies Assumption~\ref{ass},
and let $\rho_0\in L^1(\mathbb R^d)
\cap L^2(\mathbb R^d)
\cap\mathcal P_4(\mathbb R^d)$
be non-negative.
Let $\rho_\varepsilon$ be a solution to
\eqref{eq:regularized}.
Then, for any $T>0$, there exists a constant
$C_T>0$, independent of $\varepsilon$, such that
\begin{equation}
\rho_\varepsilon(t,x)\ge0
\quad\text{a.e. in }(0,T)\times\mathbb R^d,
\qquad
\int_{\mathbb R^d}\rho_\varepsilon(t,x)\,dx=1,
\end{equation}
and
\begin{equation}
\|\rho_\varepsilon\|_{
L^\infty(0,T;L^2(\mathbb R^d))}
\le C_T.
\label{L2estimate}
\end{equation}
Moreover, for every bounded open set
$\Omega\Subset\mathbb R^d$, there exists
$C_{T,\Omega}>0$, independent of $\varepsilon$,
such that
\begin{equation}
\|\partial_t\rho_\varepsilon\|_{
L^\infty(0,T;H^{-2}(\Omega))}
\le C_{T,\Omega}.
\label{time estimate}
\end{equation}
The following gradient estimates hold:
\begin{equation}
\sum_{i=1}^d
\int_0^T\int_{\mathbb R^d}
|h_i^\varepsilon(t,x)|^2
|\partial_i\rho_\varepsilon|^2\,dx\,dt
\le C_T,\quad \varepsilon^2\sum_{i=1}^d
\int_0^T\int_{\mathbb R^d}
|\partial_i\rho_\varepsilon|^2\,dx\,dt
\le C_T.
\label{local-coordinate-weight1}
\end{equation}
Finally, the second and fourth moments satisfy
\begin{align}
\operatorname*{ess\,sup}_{t\in(0,T)}
\int_{\mathbb R^d}|x|^2
\rho_\varepsilon(t,x)\,dx
&\le M_2,
\label{M2bound}\\
\operatorname*{ess\,sup}_{t\in(0,T)}
\int_{\mathbb R^d}|x|^4
\rho_\varepsilon(t,x)\,dx
&\le M_4,
\label{M4bound}
\end{align}
where $M_2$ and $M_4$ are independent
of $\varepsilon$.
\end{remark}

In what follows, we verify the spatial translation estimate and the time translation estimate separately.

\subsection{Spatial translation estimates}

To establish the spatial translation estimates, we require stronger
regularity of the solutions. We therefore derive the
following uniform $L^p(p\ge2)$ estimates for $\rho_\varepsilon$.

\begin{lemma}[$L^p$-boundedness]
\label{prop:Linfty}
Assume that $L$ satisfies Assumption~\ref{ass}, and let
$\rho_0 \in L^1(\mathbb{R}^d)\cap L^\infty(\mathbb{R}^d)$
be non-negative. Let $\rho_{\varepsilon}$ be the solution to the
regularized problem \eqref{eq:regularized}. Then, for any fixed $T>0$,
the following uniform estimates hold:
\begin{align}\label{Lpbound}
\|\rho_{\varepsilon}\|_{
L^\infty(0,T;L^p(\mathbb{R}^d))}
\le C,
\qquad p\in[2,\infty].
\end{align}
where $C$ is a constant independent of $\varepsilon$.
\end{lemma}

\begin{proof}
For any $2\le p<\infty$, testing the regularized equation \eqref{eq:regularized} with
$\rho_{\varepsilon}^{p-1}$ and integrating by parts yield
\begin{align}
\frac1p\frac{d}{dt}\|\rho_{\varepsilon}\|_{L^p(\mathbb{R}^d)}^p
&+
\frac{\sigma^2(p-1)}{2}
\sum_{i=1}^d
\int_{\mathbb{R}^d}
((h^{\varepsilon}_i)^2+\varepsilon^2)
\rho_{\varepsilon}^{p-2}|\partial_i\rho_{\varepsilon}|^2\,dx
\nonumber\\
&=
\frac{\lambda(p-1)}p
\int_{\mathbb{R}^d}
(\nabla\cdot h^{\varepsilon})\rho_{\varepsilon}^p\,dx
+
\frac{\sigma^2(p-1)}{2p}
\sum_{i=1}^d
\int_{\mathbb{R}^d}
\partial_{ii}((h^{\varepsilon}_i)^2)\rho_{\varepsilon}^p\,dx.
\label{eq:Lp-energy}
\end{align}
By the uniform bounds \eqref{h-estimate2},
there exists a constant $C_0>0$,
independent of $p$ and $\varepsilon$,
such that
\[
\frac1p\frac{d}{dt}\|\rho_{\varepsilon}\|_{L^p(\mathbb{R}^d)}^p
\le
\frac{p-1}{p}d(\lambda+\sigma^2)C_0
\|\rho_{\varepsilon}\|_{L^p(\mathbb{R}^d)}^p.
\]
Hence, by Gr\"onwall's inequality,
\[
\|\rho_{\varepsilon}(t)\|_{L^p(\mathbb{R}^d)}
\le
\exp\left(
\frac{p-1}{p}d(\lambda+\sigma^2)tC_0
\right)
\|\rho_{0,\varepsilon}\|_{L^p(\mathbb{R}^d)} \qquad \forall~0\le t\le T.
\]
Letting $p\to\infty$, we conclude that
\[
\|\rho_{\varepsilon}(t)\|_{L^\infty(\mathbb{R}^d)}
\le
e^{d(\lambda+\sigma^2)tC_0}
\|\rho_{0,\varepsilon}\|_{L^\infty(\mathbb{R}^d)}\qquad \forall~0\le t\le T.
\]
Combining the uniform $L^\infty$ bound with mass conservation,
we obtain the asserted uniform $L^p$ estimates for all
$p\in[2,\infty]$.
\end{proof}

Below is an important and useful estimate. Fix $r>0$. By the second moment estimate~\eqref{M2bound}
and Proposition~\ref{consensus point prop}, there exists
a constant $M_T>0$, independent of $\varepsilon$,
such that
\[
\sup_{0<\varepsilon\leq 1}
\sup_{t\in[0,T]}
|m_\beta^L(\rho_\varepsilon(t))|
\leq M_T.
\]

Choose $\varepsilon_r\in(0,1)$ sufficiently small so that $\varepsilon_r^{-1}>r+2+M_T.$
Then, for every $0<\varepsilon\leq\varepsilon_r$,
$x\in B_{r+2}$, and $i\in\{1,\ldots,d\}$, we have
\[
|x_i-m_{\beta,i}^L(\rho_\varepsilon(t))|
\leq r+2+M_T<\varepsilon^{-1}.
\]
Hence $\phi_{\varepsilon^{-1}}
\bigl(x_i-m_{\beta,i}^L(\rho_\varepsilon(t))\bigr)=1,$
and 
\[
h_i^\varepsilon(t,x)
=
x_i-m_{\beta,i}^L(\rho_\varepsilon(t))
\qquad\text{on }B_{r+2}.
\]
The weighted gradient estimate therefore yields
\begin{equation}
\label{local-coordinate-weight2}
\sup_{0<\varepsilon\leq\varepsilon_r}
\sum_{i=1}^d
\int_0^T\int_{B_{r+2}}
|x_i-m_{\beta,i}^L(\rho_\varepsilon(t))|^2
|\partial_i\rho_\varepsilon(t,x)|^2
\,dx\,dt
\leq C_{r,T}.
\end{equation}

%For $r>0$, let $B_r:=\{x\in\mathbb{R}^d:|x|<r\}$, by the properties of $\phi$, \eqref{local-coordinate-weight1} implies
%\begin{align}
%\sum_{i=1}^{d} \int_{0}^{T} \int_{B_r} |x_i - m_{\beta,i}^{L}(\rho_\varepsilon)|^2 |\partial_i \rho_\varepsilon|^2 \, dx \, dt \leq C.\label{local-coordinate-weight2}  
%\end{align}

Next, we employ this estimate to establish the spatial translation estimates for the regularized solutions.

\begin{lemma}[Spatial translation estimates]
\label{lem:Spatial translations estimates}
Assume that $L$ satisfies Assumption~\ref{ass}, and let the non-negative initial datum satisfy
$\rho_0 \in L^1(\mathbb{R}^d)\cap L^\infty(\mathbb{R}^d)\cap \mathcal{P}_4(\mathbb{R}^d)$. Let $\rho_{\varepsilon}$ be the solution to the
regularized problem \eqref{eq:regularized}. Then, for any $i\in\{1,\dots,d\}$ and any fixed $T, r>0$, there exists $\varepsilon_r$ such that

\begin{equation}
\sup_{0<\varepsilon\leq\varepsilon_r}
\|\rho_\varepsilon(t,x+he_i)-\rho_\varepsilon(t,x)\|_
{L^2((0,T)\times B_r)}
\longrightarrow0
\qquad\text{as }h\to0.
\label{eq:coordinate-translation}
\end{equation}
More precisely, for sufficiently small $|h|$,
\begin{equation}
\sup_{0<\varepsilon\leq\varepsilon_r}
\|\rho_\varepsilon(t,x+he_i)-\rho_\varepsilon(t,x)\|_
{L^2((0,T)\times B_r)}
\le C_{r,T,d}|h|^{\frac{1}{3}},
\label{eq:coordinate-rate}
\end{equation}
where
$C_{r,T,d}>0$ is
independent of $\varepsilon$ and $h$.
\end{lemma}

\begin{proof}
Fix $r>0$, $i\in\{1,\dots,d\}$, and let $0<|h|<\delta<1.$
For each $t\in(0,T)$, define the bad strip
\begin{equation}
S_{i,\delta}^\varepsilon(t)
:=
\left\{
x\in B_{r+1}:
|x_i-m_i^\varepsilon(t)|\le2\delta
\right\}
\label{bad-strip}
\end{equation}
and the good set
\begin{equation}
G_{i,\delta}^\varepsilon(t)
:=
\left\{
x\in B_{r+1}:
|x_i-m_i^\varepsilon(t)|>2\delta
\right\}.
\label{good-set}
\end{equation}
The strip $S_{i,\delta}^\varepsilon(t)$ is the intersection of $B_{r+1}$
with a slab of thickness $4\delta$. Therefore,
\begin{equation}
|S_{i,\delta}^\varepsilon(t)|
\le C_r\delta
\label{strip-measure}
\end{equation}
uniformly in $t$ and $\varepsilon$.

We consider
\begin{align}
\int_0^T
\int_{B_{r+1}}
|\rho_\varepsilon(t,x+he_i)-\rho_\varepsilon(t,x)|^2\,dx\,dt
=&\int_0^T
\int_{G_{i,\delta}^\varepsilon(t)}
|\rho_\varepsilon(t,x+he_i)-\rho_\varepsilon(t,x)|^2\,dx\,dt\nonumber \\
&+ \int_0^T
\int_{S_{i,\delta}^\varepsilon(t)}
|\rho_\varepsilon(t,x+he_i)-\rho_\varepsilon(t,x)|^2\,dx\,dt\label{translation second term}
\end{align}

Firstly, we estimate the translation on the good set in \eqref{good-set}. If
$x\in G_{i,\delta}^\varepsilon(t)$ and $s$ satisfy $|s|\le |h|<\delta$, then
\begin{align*}
|x_i+\operatorname{sgn}(h)s-m_i^\varepsilon(t)|
&\ge
|x_i-m_i^\varepsilon(t)|-|s|
>2\delta-\delta
=\delta.
\end{align*}
By the fundamental theorem of calculus and Cauchy--Schwarz,
\begin{align*}
|\rho_\varepsilon(t,x+he_i)-\rho_\varepsilon(t,x)|^2
=&\Big|\int_0^{|h|}
\partial_i\rho_\varepsilon(t,x+\operatorname{sgn}(h)se_i)\,ds\Big|^2
\nonumber\\
\le&
\Big(\int_0^{|h|}
|\partial_i\rho_\varepsilon(t,x+\operatorname{sgn}(h)se_i)|^2\,ds\Big)\Big(\int_0^{|h|}1^2ds\Big)
\nonumber\\
\le&
\frac{|h|}{\delta^2}
\int_0^{|h|}
|x_i+\operatorname{sgn}(h)s-m_i^\varepsilon(t)|^2
|\partial_i\rho_\varepsilon(t,x+\operatorname{sgn}(h)se_i)|^2\,ds.
\end{align*}
Integrating in $x$ over $G_{i,\delta}^\varepsilon(t)$ and then in $t$,
and using the change of variables $y=x+\operatorname{sgn}(h)se_i$, we obtain
\begin{align}
&\int_0^T
\int_{G_{i,\delta}^\varepsilon(t)}
|\rho_\varepsilon(t,x+he_i)-\rho_\varepsilon(t,x)|^2\,dx\,dt
\nonumber\\
\le&
\frac{|h|}{\delta^2}
\int_0^{|h|}
\int_0^T\!\!\int_{B_{r+2}}
|y_i-m_i^\varepsilon(t)|^2
|\partial_i\rho_\varepsilon(t,y)|^2\,dy\,dt\,ds
\le
C_{r,T,d}\frac{|h|^2}{\delta^2},
\label{eq:good-estimate}
\end{align}
where \eqref{local-coordinate-weight2} was used in the last step.

We next estimate the contribution of the bad strip in \eqref{bad-strip}. 
the uniform $L^\infty$ bound \eqref{Lpbound} and
the measure estimate \eqref{strip-measure} give
\begin{align}
&\int_0^T\!\!\int_{S_{i,\delta}^{\varepsilon}(t)}
 |\rho_\varepsilon(t,x+he_i)-\rho_\varepsilon(t,x)|^2\,dx\,dt\nonumber\\
&\qquad\le
4\|\rho_\varepsilon\|_{L^\infty((0,T)\times\mathbb R^d)}^2
\int_0^T |S_{i,\delta}^{\varepsilon}(t)|\,dt
\le C_{r,T}\delta.\label{eq:bad-estimate}
\end{align}

Substituting \eqref{eq:good-estimate}
and \eqref{eq:bad-estimate}
into \eqref{translation second term}, we obtain
\begin{equation*}
\sup_{0<\varepsilon\leq\varepsilon_r}\|\rho_\varepsilon(t,x+he_i)-\rho_\varepsilon(t,x)\|_
{L^2((0,T)\times B_r)}^2
\le
C_{r,T,d}
\left(
\frac{|h|^2}{\delta^2}
+
\delta
\right).
\end{equation*}
We optimize the right-hand side by choosing $\delta
=
|h|^{\frac{2}{3}}.$
Then,
\[
\sup_{0<\varepsilon\leq\varepsilon_r}\|\rho_\varepsilon(t,x+he_i)-\rho_\varepsilon(t,x)\|_
{L^2((0,T)\times B_r)}
\le
C_{r,T,d}
|h|^{\frac{1}{3}},
\]
uniformly for \(0<\varepsilon\le\varepsilon_r\), proving the assertion.
\end{proof}

\begin{remark}
The spatial translation estimates established above
hold uniformly for sufficiently small $\varepsilon$,
where the upper bound $\varepsilon_r$ may depend
on the radius $r$ of the spatial domain.
This dependence does not affect the subsequent
compactness argument. Indeed, for any sequence
$\varepsilon_k\to0$ and any fixed $r>0$, there
exists $k_r\in\mathbb{N}$ such that
$\varepsilon_k\leq\varepsilon_r$ for all
$k\geq k_r$.
Consequently, the spatial translation estimates
hold uniformly along the tail of the sequence
on every fixed bounded domain, which is sufficient
for the application of the
Fr\'echet--Kolmogorov compactness criterion.
\end{remark}

In the following, we extend the coordinate-wise translation estimates to arbitrary spatial translations by decomposing a general translation vector into its coordinate components.

\begin{lemma}[General spatial translations]
\label{lem:general-spatial}
Assume that $L$ satisfies Assumption~\ref{ass},
and let
$\rho_0 \in L^1(\mathbb{R}^d)
\cap L^\infty(\mathbb{R}^d)
\cap \mathcal{P}_4(\mathbb{R}^d)$
be non-negative.
Let $\rho_{\varepsilon}$ be the solution to
the regularized problem \eqref{eq:regularized}.
Then, for any fixed $T,r>0$,
there exists $\varepsilon_r>0$ such that,
for all sufficiently small $|y|$,
\begin{equation}
\sup_{0<\varepsilon\leq\varepsilon_r}
\|\rho_\varepsilon(t,x+y)
-\rho_\varepsilon(t,x)\|_{
L^2((0,T)\times B_r)}
\le
C_{r,T,d}
|y|^{\frac{1}{3}}.
\label{eq:general-spatial}
\end{equation}
\end{lemma}
\begin{proof}
Choose $\varepsilon_r>0$ sufficiently small so that
Lemma~\ref{lem:Spatial translations estimates} applies on $B_{r+1}$ for every coordinate direction.
Write
\[
y=(y_1,\dots,y_d)=y_1e_1+\cdots+y_de_d
\]
and decompose the translation into coordinate
translations:
\begin{align*}
\rho_\varepsilon(t,x+y)-\rho_\varepsilon(t,x)
&=
\sum_{i=1}^d
\Big[
\rho_\varepsilon
\bigl(t,x+y_1e_1+\cdots+y_ie_i\bigr)
\\
&\qquad\qquad-
\rho_\varepsilon
\bigl(t,x+y_1e_1+\cdots+y_{i-1}e_{i-1}\bigr)
\Big].
\end{align*}
For $|y|<1$, all intermediate points
$x+\sum_{j=1}^i y_je_j$, with $x\in B_r$,
belong to $B_{r+1}$.
Thus, by the triangle inequality, a change
of variables, and
Lemma~\ref{lem:Spatial translations estimates},
we obtain, for all sufficiently small $|y|$,
\begin{align*}
&\|\rho_\varepsilon(\cdot,\cdot+y)
-\rho_\varepsilon\|_{L^2((0,T)\times B_r)}
\le
C_{r,T}\sum_{i=1}^d|y_i|^{\frac{1}{3}}
\le
C_{r,T,d}|y|^{\frac{1}{3}},
\end{align*}
uniformly for
$0<\varepsilon\le\varepsilon_r$. Taking the supremum over
$0<\varepsilon\le\varepsilon_r$
proves \eqref{eq:general-spatial}.
\end{proof}

\subsection{Temporal translation estimate}

We now investigate the temporal regularity of the approximate solutions. The following lemma provides the corresponding uniform temporal translation estimate.

%The preceding estimate controls time translations only in a negative Sobolev space. We combine it with the spatial translation estimate to obtain the $L^2$ time equicontinuity required by Fr\'echet--Kolmogorov.

\begin{lemma}[Uniform temporal translations]\label{lem:temporal-translations}
Assume that $L$ satisfies Assumption \ref{ass},
and let
$\rho_0\in L^1(\mathbb R^d)
\cap L^\infty(\mathbb R^d)
\cap\mathcal P_4(\mathbb R^d)$
be non-negative.
Let $\rho_\varepsilon$ be the solution
to the regularized problem \eqref{eq:regularized}.
Then, for any $T,r>0$, there exists
$\varepsilon_r>0$ such that
\begin{equation}
\lim_{\tau\to0}
\sup_{0<\varepsilon\leq\varepsilon_r}
\|\rho_\varepsilon(t+\tau,x)
-\rho_\varepsilon(t,x)\|
_{L^2(I_\tau\times B_r)}
=0,
\label{3.22}
\end{equation}
where, for $|\tau|<T$,
\[
I_\tau=
\begin{cases}
(0,T-\tau), & 0\leq\tau<T,\\
(-\tau,T), & -T<\tau<0.
\end{cases}
\]
\end{lemma}

%\begin{lemma}[Uniform temporal translations]
%\label{lem:temporal-translations}
%Assume that $L$ satisfies Assumption~\ref{ass}, and let$\rho_0 \in L^1(\mathbb{R}^d)\cap L^\infty(\mathbb{R}^d)\cap \mathcal{P}_4(\mathbb{R}^d)$be non-negative. Let $\rho_{\varepsilon}$ be the solution to the regularized problem \eqref{eq:regularized}. Then, for any $T,r>0$ and $0<\tau<T$, there exists $\varepsilon_r$ such that
%\begin{equation}\lim_{\tau\to0^+}\sup_{0<\varepsilon<\varepsilon_r}\|\rho_\varepsilon(t+\tau,x)-\rho_\varepsilon(t,x)\|_{L^2((0,T-\tau)\times B_r)}=0.\label{eq:temporal-translations}\end{equation}\end{lemma}

\begin{proof}
Fix \(T,r>0\), and let \(\varepsilon_r>0\) be the threshold given by Lemma~\ref{lem:general-spatial} on \(B_r\). Choose \(\chi\in C_c^\infty(B_{r+2})\) such that \(0\le\chi\le1\) and \(\chi\equiv1\) on \(B_{r+1}\). For \(0<\eta<1/2\), set \(J_\eta(x)=\eta^{-d}J(x/\eta)\), where \(J\) is the mollifier introduced in \eqref{J}, and define \(\rho_\varepsilon^\eta=J_\eta*(\chi\rho_\varepsilon)\).

We first consider $0<\tau<T$. By the triangle inequality, we have
\begin{align}\label{time triangle inequality}
\|\rho_\varepsilon(t+\tau,x)-\rho_\varepsilon(t,x)\|_
{L^2((0,T-\tau)\times B_r)}
\le&
\|\rho_\varepsilon(t+\tau,x)
-\rho_\varepsilon^\eta(t+\tau,x)\|_
{L^2((0,T-\tau)\times B_r)}
\nonumber\\
&+
\|\rho_\varepsilon^\eta(t+\tau,x)
-\rho_\varepsilon^\eta(t,x)\|_
{L^2((0,T-\tau)\times B_r)}
\nonumber\\
&+
\|\rho_\varepsilon^\eta(t,x)-\rho_\varepsilon(t,x)\|_
{L^2((0,T-\tau)\times B_r)}.
\end{align}

As for the first and third terms, since $\chi\equiv1$ on $B_{r+1}$ and
$\operatorname{supp}J_\eta\subset B_\eta$, for $x\in B_r$
we have
\[
\rho_\varepsilon^\eta(\cdot,x)-\rho_\varepsilon(\cdot,x)
=
\int_{\mathbb{R}^d}
J_\eta(y)
\bigl[
\chi\rho_\varepsilon(\cdot,x-y)-\rho_\varepsilon(\cdot,x)
\bigr]\,dy .
\]
By Minkowski's inequality and the uniform spatial
translation estimate of Lemma~\ref{lem:general-spatial},
\begin{align}
\|\rho_\varepsilon^\eta-\rho_\varepsilon\|_
{L^2((0,T)\times B_r)}
&\le
\int_{\mathbb{R}^d}
J_\eta(y)
\|\rho_\varepsilon(\cdot,x-y)-\rho_\varepsilon(\cdot,x)\|_
{L^2((0,T)\times B_r)}
\,dy
\nonumber\\
&\le \sup_{|y|\le \eta}\|\rho_\varepsilon(\cdot,x-y)-\rho_\varepsilon(\cdot,x)\|_
{L^2((0,T)\times B_r)}
\nonumber\\
&\le C_{r,T,d}\eta^{\frac{1}{3}}
\longrightarrow0,
\quad\text{as }\eta\to0,
\label{eq:mollification-error}
\end{align}

We next estimate the second term of \eqref{time triangle inequality}. For $0<\tau<T$,
\[
\rho_\varepsilon^\eta(t+\tau,x)-\rho_\varepsilon^\eta(t,x)
=
J_\eta*
\left[
\chi\bigl(
\rho_\varepsilon(t+\tau,x)-\rho_\varepsilon(t,x)
\bigr)
\right].
\]
For each fixed $\eta>0$, the operator
\[
v\longmapsto J_\eta*(\chi v)
\]
is continuous from $H^{-2}(B_{r+2})$ into $L^2(B_r)$.
Indeed, for $x\in B_r$,
\begin{align*}
|J_\eta*(\chi v)(x)|
&=|\left\langle
v,\,
\chi(\cdot)J_\eta(x-\cdot)
\right\rangle|\\ 
&\le
\|v\|_{H^{-2}(B_{r+2})}
\,
\|\chi(\cdot)J_\eta(x-\cdot)\|_{H^2(B_{r+2})},  
\end{align*}
Integrating the squared pointwise estimate over $x\in B_r$
and taking the square root, we obtain
\begin{equation}\label{eq:mollifier-Hminus2}
\|J_\eta*(\chi v)\|_{L^2(B_r)}
\le
C_r\eta^{-d/2-2}\|v\|_{H^{-2}(B_{r+2})}.
\end{equation}
Applying this 
gives
\begin{align*}
\|\rho_\varepsilon^\eta(t+\tau,x)
-\rho_\varepsilon^\eta(t,x)\|_{L^2(B_r)}
&\le
C_{r}\eta^{-d/2-2}
\|\rho_\varepsilon(t+\tau,x)-\rho_\varepsilon(t,x)\|_
{H^{-2}(B_{r+2})}
\nonumber\\
&\le
C_{r}\eta^{-d/2-2}\int_t^{t+\tau}
\|\partial_s\rho_\varepsilon(s,x)\|_
{H^{-2}(B_{r+2})}\,ds
\\
&\le C_{r}\eta^{-d/2-2}\tau.
\label{eq:mollified-time-translation}
\end{align*}
\eqref{time estimate} is used here. Integrating with respect to $t\in(0,T-\tau)$ yields
\begin{equation}
\|\rho_\varepsilon^\eta(t+\tau,x)
-\rho_\varepsilon^\eta(t,x)\|_
{L^2((0,T-\tau)\times B_r)}
\le
C_{r,T}\eta^{-d/2-2}\,\tau .
\label{eq:mollified-time-L2}
\end{equation}

Finally, substituting \eqref{eq:mollification-error}
and \eqref{eq:mollified-time-L2} into
\eqref{time triangle inequality}, we obtain
\begin{align*}
&\|\rho_\varepsilon(t+\tau,x)
-\rho_\varepsilon(t,x)\|_
{L^2((0,T-\tau)\times B_r)}
\le
C_{r,T,d}
\left(
\eta^{\frac{1}{3}}
+\eta^{-\frac d2-2}\tau
\right).
\end{align*}
Choosing $\eta
=
\tau^{
\frac{1}{
{\frac{1}{3}}+\frac d2+2
}},$
for sufficiently small $\tau>0$, we obtain
\begin{align*}
&\sup_{0<\varepsilon\le\varepsilon_r}
\|\rho_\varepsilon(t+\tau,x)
-\rho_\varepsilon(t,x)\|_
{L^2((0,T-\tau)\times B_r)}
\le
C_{r,T,d}
\tau^{
\frac{2}{14+3d}}
\longrightarrow0
\qquad\text{as }\tau\to0^+.
\end{align*}

The same estimate holds for negative time
translations by a change of variables.
More precisely, for $-T<\tau<0$,
setting $s=t+\tau$, we have
\[
\begin{aligned}
&\|\rho_\varepsilon(t+\tau)-\rho_\varepsilon(t)\|
_{L^2((-\tau,T)\times B_r)}=
\|\rho_\varepsilon(s+|\tau|)-\rho_\varepsilon(s)\|
_{L^2((0,T-|\tau|)\times B_r)}.
\end{aligned}
\]

Consequently, for both positive and
negative time translations,
\[
\lim_{\tau\to0}
\sup_{0<\varepsilon\leq\varepsilon_r}
\|\rho_\varepsilon(t+\tau,x)-\rho_\varepsilon(t,x)\|
_{L^2(I_\tau\times B_r)}
=0,
\]
where
\[
I_\tau=
\begin{cases}
(0,T-\tau), & 0\leq\tau<T,\\
(-\tau,T), & -T<\tau<0.
\end{cases}
\]
The proof is complete.

\end{proof}

\section{Compactness and existence} \label{existence}

In this section, we establish strong compactness in the whole space
and obtain global weak solutions. We first derive a uniform
$L^2$-tail estimate, which complements the translation estimates
in the application of the Fr\'echet--Kolmogorov compactness criterion.
\begin{lemma}[Uniform $L^2$-tail estimate]
\label{lem:L2_tail}
Assume that $L$ satisfies Assumption~\ref{ass},
and let the non-negative initial datum satisfy $\rho_0 \in L^1(\mathbb{R}^d)
\cap L^\infty(\mathbb{R}^d)
\cap \mathcal{P}_4(\mathbb{R}^d).$
Let $\rho_\varepsilon$ be the solution to the
regularized problem. Then, for any fixed $T>0$,
we have
\begin{equation}
\lim_{r\to\infty}
\sup_{0<\varepsilon\le1}
\int_0^T\int_{|x|>r}
|\rho_\varepsilon(t,x)|^2\,dx\,dt=0.
\label{eq:L2_tail}
\end{equation}
\end{lemma}

\begin{proof}
For any $r>0$, using the non-negativity
of $\rho_\varepsilon$ and the uniform
$L^\infty$ estimate \eqref{Lpbound}, we obtain
\begin{align*}
\int_0^T\int_{|x|>r}
|\rho_\varepsilon(t,x)|^2\,dx\,dt
&\le
\|\rho_\varepsilon\|_{
L^\infty((0,T)\times\mathbb{R}^d)}
\int_0^T\int_{|x|>r}
\rho_\varepsilon(t,x)\,dx\,dt\\
&\le
\frac{C}{r^4}
\int_0^T\int_{\mathbb{R}^d}
|x|^4\rho_\varepsilon(t,x)\,dx\,dt\le
\frac{C T M_4}{r^4}.
\end{align*}
Here, the last inequality follows from
the uniform fourth moment bound \eqref{M4bound}. Since the right-hand side is independent
of $\varepsilon$, taking the supremum over
$0<\varepsilon\le1$ yields
\[
\sup_{0<\varepsilon\le1}
\int_0^T\int_{|x|>r}
|\rho_\varepsilon(t,x)|^2\,dx\,dt
\le \frac{CT M_4}{r^4}.
\]
Letting $r\to\infty$, we conclude that
\[
\lim_{r\to\infty}
\sup_{0<\varepsilon\le1}
\int_0^T\int_{|x|>r}
|\rho_\varepsilon(t,x)|^2\,dx\,dt=0.
\]
This completes the proof.
\end{proof}

Combining the uniform $L^2$-tail estimate with the
spatial and temporal translation estimates established
in the previous section, we obtain the following
strong compactness result on the whole space.

\begin{proposition}[Strong compactness on the whole space]
Assume that $L$ satisfies Assumption~\ref{ass},
and let the non-negative initial datum satisfy
$\rho_0\in L^1(\mathbb R^d)
\cap L^\infty(\mathbb R^d)
\cap\mathcal P_4(\mathbb R^d)$.
Let $\rho_\varepsilon$ be the solution
to the regularized problem
\eqref{eq:regularized}.
Let $\{\varepsilon_n\}\subset(0,1]$ be any sequence
satisfying $\varepsilon_n\to0$ as $n\to\infty$.
Then, for any fixed $T>0$, there exist
a subsequence, still denoted by
$\{\varepsilon_n\}$ and a non-negative function
$\rho\in L^2((0,T)\times\mathbb R^d)$
such that
\[
\rho_{\varepsilon_n}\to\rho
\quad\text{strongly in }
L^2((0,T)\times\mathbb R^d).
\]
\end{proposition}

%\begin{proposition}
%\label{prop:strong_compactness}
%Assume that $L$ satisfies Assumption~\ref{ass},
%and let the non-negative initial datum satisfy $\rho_0\in L^1(\mathbb{R}^d)
%\cap L^\infty(\mathbb{R}^d)
%\cap\mathcal{P}_4(\mathbb{R}^d).$
%Let $\rho_\varepsilon$ be the solution to the regularized problem~\eqref{eq:regularized}.
%Then, for any fixed $T>0$, the family$\{\rho_\varepsilon\}_{\varepsilon>0}$is relatively compact in$L^2((0,T)\times\mathbb{R}^d)$. Consequently, there exist a subsequence$\varepsilon_k\to0$ and a non-negative function$\rho\in L^2((0,T)\times\mathbb{R}^d)$such that\begin{equation}\rho_{\varepsilon_k}\to\rho\quad\text{strongly in }L^2((0,T)\times\mathbb{R}^d).\label{eq:strong_L2}\end{equation}\end{proposition}

\begin{proof}
Fix $T>0$ and let $\{\varepsilon_n\}$ be any
sequence satisfying $\varepsilon_n\to0$.
We apply the Fr\'echet--Kolmogorov compactness
criterion in Theorem~\ref{thm:FK} with $p=2$
to the sequence $\{\rho_{\varepsilon_n}\}_{n\ge1}$.

First, by the uniform $L^2$ estimate
\eqref{L2estimate}, we have
\[
\sup_{n\ge1}
\|\rho_{\varepsilon_n}\|_{L^2((0,T)\times\mathbb R^d)}
\le C_T.
\]
Thus, the sequence is uniformly bounded
in $L^2((0,T)\times\mathbb R^d)$.

Next, we verify condition (i).
Fix $r>0$.
By Lemmas~\ref{lem:general-spatial}
and \ref{lem:temporal-translations},
there exists $\varepsilon_r>0$ such that
\[
\lim_{|\xi|+|\tau|\to0}
\sup_{0<\varepsilon\le\varepsilon_r}
\left(
\|\rho_\varepsilon(t,x+\xi)
-\rho_\varepsilon(t,x)\|_
{L^2((0,T)\times B_r)}
+
\|\rho_\varepsilon(t+\tau,x)
-\rho_\varepsilon(t,x)\|_
{L^2(I_\tau\times B_r)}
\right)=0.
\]
Since $\varepsilon_n\to0$, there exists
$N_r\in\mathbb N$ such that
$\varepsilon_n\le\varepsilon_r$
for all $n\ge N_r$.

On the other hand, the remaining finitely
many functions
$\rho_{\varepsilon_1},\ldots,
\rho_{\varepsilon_{N_r-1}}$
belong to $L^2((0,T)\times\mathbb R^d)$.
By the continuity of translations in $L^2$,
their spatial and temporal translation
differences also converge to zero
as $|\xi|+|\tau|\to0$.
Consequently,
\[
\begin{aligned}
\lim_{|\xi|+|\tau|\to0}
\sup_{n\ge1}
\Big(
&\|\rho_{\varepsilon_n}(t,x+\xi)
-\rho_{\varepsilon_n}(t,x)\|_
{L^2((0,T)\times B_r)}\\
+&
\|\rho_{\varepsilon_n}(t+\tau,x)
-\rho_{\varepsilon_n}(t,x)\|_
{L^2(I_\tau\times B_r)}
\Big)=0.
\end{aligned}
\]
Therefore, condition (i) holds.

We now verify condition (ii).
By the uniform $L^2$-tail estimate
in Lemma~\ref{lem:L2_tail}, we have
\[
\lim_{r\to\infty}
\sup_{n\ge1}
\int_0^T\int_{|x|>r}
|\rho_{\varepsilon_n}(t,x)|^2
\,dx\,dt=0.
\]
Hence, for any $\eta>0$, there exists
$r>0$ such that
\[
\sup_{n\ge1}
\|\rho_{\varepsilon_n}\|_
{L^2((0,T)\times
(\mathbb R^d\setminus B_r))}
<\eta.
\]
Thus, condition (ii) is satisfied.

Consequently, all the assumptions of
Theorem~\ref{thm:FK} are satisfied.
Therefore, the sequence
$\{\rho_{\varepsilon_n}\}_{n\ge1}$
is relatively compact in $L^2((0,T)\times\mathbb R^d)$.
Hence, there exist a subsequence,
still denoted by $\{\varepsilon_n\}$,
and a function $\rho\in L^2((0,T)\times\mathbb R^d)$
such that
\[
\rho_{\varepsilon_n}\to\rho
\quad\text{strongly in }L^2((0,T)\times\mathbb R^d).
\]
Since $\rho_{\varepsilon_n}\ge0$,
we also have $\rho\ge0$ almost everywhere.
This completes the proof.
\end{proof}

\begin{remark}
The preceding proposition highlights a key difference
between the present analysis and the isotropic setting.
Here, the Fr\'echet--Kolmogorov compactness criterion
replaces the Aubin--Lions compactness argument, without
requiring uniform local $W^{1,\gamma}$ estimates.
Instead, strong compactness is obtained directly from
the spatial and temporal translation estimates, together
with the uniform $L^2$-tail control.
\end{remark}

%We next extend this local convergence to the whole space. In Reference ~\cite{XXX}, the local compactness of the isotropic CBO equation was obtained by means of a local Sobolev estimate and the Aubin-Lions lemma, as in $\varepsilon\to 0$, and the following convergence holds
%\begin{align}\label{Rtovar}\rho_{\varepsilon} \to \rho\qquad \mbox{strongly in} \ L^2(0,T; L^q(B_r)), ~~q\in[1, \frac{d\gamma}{d-\gamma}).\end{align}
%where $\gamma\in(1,\frac{2d}{2+d})$. Since, in the proof of global strong convergence in \cite[Lemma~3.1]{XXX}, the property of local strong convergence is used only to obtain almost convergence everywhere (up to a subsequence), our local
%$L^2$-convergence 
%\[
%\rho_{\varepsilon_k}\to\rho\qquad\text{strongly in }L^2((0,T)\times B_r).\]
%provides exactly the same property. Hence, the remaining argument of \cite[Lemma~3.1]{XXX} applies without modification.

We next collect the convergence properties of the
approximate solutions and the estimates inherited
by their limit, which will be used to pass to the
limit in the regularized equation.

\begin{corollary}[Convergence properties of the approximate solutions]
\label{Convergence properties of the approximate solutions}
Assume that $L$ satisfies Assumption~\ref{ass},
and let the non-negative initial datum satisfy
$\rho_0\in L^1(\R^d)\cap L^\infty(\R^d)
\cap\mathcal{P}_4(\R^d)$.
Let $\rho_\varepsilon$ be the solution to the
regularized problem \eqref{eq:regularized}.
Then, for any fixed $T>0$, there exist a sequence
$\varepsilon_k\to0$ and a non-negative function $\rho\in L^\infty(0,T;L^2(\R^d))
\cap L^\infty((0,T)\times\R^d)$
such that
\begin{align}
\rho_{\varepsilon_k}
&\to\rho
&&\textup{strongly in }
L^2((0,T)\times\R^d),
\label{rhouniform}\\
\rho_{\varepsilon_k}
&\overset{*}{\rightharpoonup}\rho
&&\textup{in }
L^\infty(0,T;L^p(\R^d)),
\quad\text{for every }1<p<\infty,
\label{Rntovar3}\\
\partial_t\rho_{\varepsilon_k}
&\overset{*}{\rightharpoonup}\partial_t\rho
&&\textup{in }
L^\infty(0,T;H^{-2}(\Omega)),
\quad\text{for every bounded open }
\Omega\Subset\R^d,
\label{Rntovar4}
\end{align}
\end{corollary}

Moreover, as in Lemma 3.2 of \cite{huang2026mathematical}, we have the following conclusion.

\begin{lemma}\label{W2estimate}
Assume that $L$ satisfies Assumption~\ref{ass},
and let the non-negative initial datum satisfy
$\rho_0\in L^1(\mathbb R^d)
\cap L^\infty(\mathbb R^d)
\cap\mathcal P_4(\mathbb R^d)$.
Let $\rho_{\varepsilon_k}$ and $\rho$
be as in Corollary~\ref{Convergence properties of the approximate solutions}.
Then, after extracting a further subsequence
if necessary, we have
\begin{align}
\rho_{\varepsilon_k}
&\longrightarrow\rho
\quad\text{strongly in }
L^1((0,T)\times\mathbb R^d),
\label{eq:globalL1}\\
W_2(\rho_{\varepsilon_k}(t),\rho(t))
&\longrightarrow0
\quad\text{for almost every }t\in(0,T).
\end{align}
In particular, by the $W_2$-Lipschitz estimate in Proposition~\ref{consensus point prop}
and the uniform second-moment bounds, we obtain
\begin{align}\label{consensus point estimate}
&|m_\beta^L(\rho_{\varepsilon_k}(t))
-m_\beta^L(\rho(t))|\le C_m
W_2(\rho_{\varepsilon_k}(t),\rho(t))
\longrightarrow0
\end{align}
for almost every $t\in(0,T)$.
\end{lemma}

\begin{proof}
By \eqref{rhouniform}, after passing to a further subsequence, still denoted
by $\{\varepsilon_k\}$, we have
\[
\|\rho_{\varepsilon_k}(t)-\rho(t)\|_{L^2(\mathbb R^d)}
\longrightarrow 0
\quad\text{for almost every }t\in(0,T),
\]
and
\[
\rho_{\varepsilon_k}(t,x)\longrightarrow\rho(t,x)
\quad\text{for almost every }(t,x)\in(0,T)\times\mathbb R^d.
\]
By Fatou's lemma and the uniform moment bounds
\eqref{M2bound} and \eqref{M4bound},
we obtain
\begin{align}
\int_{\mathbb R^d}|x|^2\rho(t,x)\,dx
&\le
\liminf_{k\to\infty}
\int_{\mathbb R^d}|x|^2
\rho_{\varepsilon_{k}}(t,x)\,dx
\le M_2,
\label{M2bound-limit}\\
\int_{\mathbb R^d}|x|^4\rho(t,x)\,dx
&\le
\liminf_{k\to\infty}
\int_{\mathbb R^d}|x|^4
\rho_{\varepsilon_{k}}(t,x)\,dx
\le M_4.
\label{M4bound-limit}
\end{align}
These estimates hold for almost every
$t\in(0,T)$.

We next prove the strong $L^1$ convergence.
First, for any $r>0$, we decompose
\[
\begin{aligned}
\|\rho_{\varepsilon_k}(t)-\rho(t)\|_
{L^1(\mathbb R^d)}
&=
\int_{B_r}
|\rho_{\varepsilon_k}(t,x)-\rho(t,x)|\,dx
+
\int_{\mathbb R^d\setminus B_r}
|\rho_{\varepsilon_k}(t,x)-\rho(t,x)|\,dx.
\end{aligned}
\]
For the first term, H\"older's inequality yields
\[
\int_{B_r}
|\rho_{\varepsilon_k}(t,x)-\rho(t,x)|\,dx
\le |B_r|^{1/2}
\|\rho_{\varepsilon_k}(t)-\rho(t)\|_{L^2(\mathbb R^d)}
\longrightarrow 0
\quad\text{as }k\to\infty.
\]
For the second term, using the
non-negativity of the densities and the fourth moment estimates \eqref{M4bound} and \eqref{M4bound-limit}, we obtain
\[
\begin{aligned}
\int_{\mathbb R^d\setminus B_r}
|\rho_{\varepsilon_k}(t,x)-\rho(t,x)|\,dx&\leq
\int_{\mathbb R^d\setminus B_r}\rho_{\varepsilon_k}(t,x)\,dx
+
\int_{\mathbb R^d\setminus B_r}\rho(t,x)\,dx
\\
&\leq
\frac{1}{r^4}
\int_{\mathbb R^d\setminus B_r}|x|^4
\rho_{\varepsilon_k}(t,x)\,dx + 
\frac{1}{r^4}
\int_{\mathbb R^d\setminus B_r}|x|^4
\rho(t,x)\,dx.
\leq\frac{2M_4}{r^4}.
\end{aligned}
\]
Thus,
\begin{equation}
\begin{aligned}
&\|\rho_{\varepsilon_k}(t)-\rho(t)\|_
{L^1(\mathbb R^d)}
\leq
|B_r|^{1/2}
\|\rho_{\varepsilon_k}(t)-\rho(t)\|_
{L^2(\mathbb R^d)}
+\frac{2M_4}{r^4}.
\end{aligned}
\label{eq:L1-estimate}
\end{equation}
For almost every $t\in(0,T)$, taking the upper limit as
$k\to\infty$ in \eqref{eq:L1-estimate} gives, for every $r>0$,
\[
\limsup_{k\to\infty}
\|\rho_{\varepsilon_k}(t)-\rho(t)\|_{L^1(\mathbb R^d)}
\le \frac{2M_4}{r^4}.
\]
Letting $r\to\infty$ yields (4.11).
\begin{equation}
\rho_{\varepsilon_k}(t)\longrightarrow\rho(t)
\quad\text{strongly in }L^1(\mathbb R^d).
\label{eq:global-L1}
\end{equation}
Consequently,
\[
\left|\int_{\mathbb R^d}\rho(t,x)\,dx-\int_{\mathbb R^d}\rho_{\varepsilon_k}(t,x)\,dx\right|
\le
\|\rho(t)-\rho_{\varepsilon_k}(t)\|_{L^1(\mathbb R^d)}
\longrightarrow0.
\]
Thus $\rho(t,x)\,dx$ is a probability measure
for almost every $t\in(0,T)$.

Moreover, since both $\rho_{\varepsilon_k}(t)$
and $\rho(t)$ are non-negative probability densities,
we have
\[
\|\rho_{\varepsilon_k}(t)-\rho(t)\|_{L^1(\mathbb R^d)}
\leq
\|\rho_{\varepsilon_k}(t)\|_{L^1(\mathbb R^d)}
+
\|\rho(t)\|_{L^1(\mathbb R^d)}\le C
\]
for almost every $t\in(0,T)$.
By the dominated convergence theorem,
\[
\begin{aligned}
&\|\rho_{\varepsilon_k}-\rho\|_
{L^1((0,T)\times\mathbb R^d)}=
\int_0^T
\|\rho_{\varepsilon_k}(t)-\rho(t)\|_
{L^1(\mathbb R^d)}\,dt
\longrightarrow0.
\end{aligned}
\]

Furthermore, for every bounded continuous
function $\varphi$ on $\mathbb R^d$,
\[
\begin{aligned}
&\left|
\int_{\mathbb R^d}\varphi(x)
\bigl(\rho_{\varepsilon_k}(t,x)-\rho(t,x)\bigr)
\,dx\right|\leq
\|\varphi\|_{L^\infty(\mathbb R^d)}
\|\rho_{\varepsilon_k}(t)-\rho(t)\|_
{L^1(\mathbb R^d)}
\longrightarrow0.
\end{aligned}
\]
Therefore,
\begin{equation}
\rho_{\varepsilon_k}(t,x)\,dx
\rightharpoonup
\rho(t,x)\,dx
\quad\text{narrowly}.
\label{eq:narrow-convergence}
\end{equation}

Next, we establish the convergence of the second moments.
\[
\int_{\mathbb R^d}|x|^2
\rho_{\varepsilon_k}(t,x)\,dx
\longrightarrow
\int_{\mathbb R^d}|x|^2
\rho(t,x)\,dx.
\]
For any $r>0$, we write
\begin{align}\label{220}
\left|
\int_{\mathbb R^d}|x|^2
\bigl(\rho_{\varepsilon_k}(t,x)-\rho(t,x)\bigr)
\,dx\right|
&\leq
\int_{B_r}|x|^2
|\rho_{\varepsilon_k}(t,x)-\rho(t,x)|\,dx
\nonumber\\
&+
\int_{\mathbb R^d\setminus B_r}|x|^2
\rho_{\varepsilon_k}(t,x)\,dx
+
\int_{\mathbb R^d\setminus B_r}|x|^2\rho(t,x)\,dx.
\end{align}
On the bounded domain $B_r$, we have
$|x|^2\leq r^2$. Hence,

\begin{align}\label{221}
&\int_{B_r}|x|^2
|\rho_{\varepsilon_k}(t,x)-\rho(t,x)|\,dx\leq
r^2\|\rho_{\varepsilon_k}(t)-\rho(t)\|_
{L^1(\mathbb R^d)}
\longrightarrow0.
\end{align}
On the other hand,
\begin{align}\label{222}
\int_{\mathbb R^d\setminus B_r}|x|^2
\rho_{\varepsilon_k}(t,x)\,dx
&\leq
\frac1{r^2}
\int_{\mathbb R^d\setminus B_r}|x|^4
\rho_{\varepsilon_k}(t,x)\,dx
\leq\frac{M_4}{r^2}.
\end{align}
Similarly,
\begin{align}\label{223}
\int_{\mathbb R^d\setminus B_r}|x|^2\rho(t,x)\,dx
\leq\frac{M_4}{r^2}. 
\end{align}
Combining \eqref{220}--\eqref{223} yields
\[
\begin{aligned}
&\left|
\int_{\mathbb R^d}|x|^2
\rho_{\varepsilon_k}(t,x)\,dx
-
\int_{\mathbb R^d}|x|^2
\rho(t,x)\,dx
\right|
\leq
r^2\|\rho_{\varepsilon_k}(t)-\rho(t)\|_
{L^1(\mathbb R^d)}
+\frac{2M_4}{r^2}.
\end{aligned}
\]
First letting $k\to\infty$ and then
$r\to\infty$, we conclude that
\begin{equation}
\int_{\mathbb R^d}|x|^2
\rho_{\varepsilon_k}(t,x)\,dx
\longrightarrow
\int_{\mathbb R^d}|x|^2
\rho(t,x)\,dx.
\label{eq:second moment-convergence}
\end{equation}

Combining \eqref{eq:narrow-convergence}
and \eqref{eq:second moment-convergence},
the characterization of convergence in
$\mathcal P_2(\mathbb R^d)$ yields
\[
W_2(\rho_{\varepsilon_k}(t),\rho(t))
\longrightarrow0
\quad\text{as }k\to\infty.
\]
%Consequently, by \eqref{eq:narrow-convergence}, the probability measures $\rho_{\varepsilon_k}(t,x)\,dx$
%converge narrowly to $\rho(t,x)\,dx$.
%Moreover, by \eqref{eq:second moment-convergence}, their second moments converge.
%It is well known that convergence in $\mathcal P_2(\mathbb R^d)$ with respect to the 2-Wasserstein distance is equivalent to narrow convergence together with convergence of the second moments. Therefore,
%\[W_2(\rho_{\varepsilon_k}(t),\rho(t))\longrightarrow0\quad\text{for almost every }t\in(0,T) \quad\text{as }k\to\infty.\]
Finally, by the stability estimate for the consensus point $m_\beta^L$ in Proposition~\ref{consensus point prop}, as $k\to\infty$, 
\begin{align*}
|m_\beta^L(\rho_{\varepsilon_k})- m_\beta^L(\rho)|\le C_m W_2(\rho_{\varepsilon_k}(t),\rho(t))
\longrightarrow0,
\end{align*}
for almost every $t\in(0,T)$. This completes the proof.
\end{proof}

With these convergence properties at hand, we now
pass to the limit in the weak formulation of the
regularized problem and show that the limit function
is a weak solution of the original equation.

\noindent\textbf{Proof of Theorem~\ref{thm:main}.} 
Throughout the proof, $\varepsilon\to0$ along the subsequence selected in Lemma~\ref{W2estimate}.
Take any test function $\varphi\in C_c^\infty([0,T)\times\R^d)$. Choose a bounded domain
$\Omega\subset\R^d$ such that
$
\operatorname{supp}\varphi
\subset [0,T)\times\Omega.
$
%Taking any test function  $\varphi\in C_c^\infty([0,T)\times\R^d)$ and $\Omega:=\supp\varphi$.
The weak formulation of the regularized equation \eqref{eq:regularized} reads as follows
\begin{align}
-\int_{\mathbb{R}^d}
\rho_{0,\varepsilon}(x)\varphi(0,x)\,dx ={}&
\int_0^T\int_{\mathbb{R}^d}
\rho_\varepsilon\,\partial_t\varphi\,dx\,dt
-\lambda\sum_{i=1}^d
\int_0^T\int_{\mathbb{R}^d}
h_i^\varepsilon(t,x)\rho_\varepsilon\,\partial_i\varphi\,dx\,dt
\nonumber\\
&+\frac{\sigma^2}{2}\sum_{i=1}^d
\int_0^T\int_{\mathbb{R}^d}
\bigl(h_i^\varepsilon(t,x)^2+\varepsilon^2\bigr)
\rho_\varepsilon\,\partial_{ii}\varphi\,dx\,dt.
\label{eq:weak-regularized}
\end{align}

Recall that for each 
$i=1,\ldots,d$, the regularized coefficient is defined as
\[
h_i^{\varepsilon}(t,x)
=
\left(x_i-m_{\beta,i}^L
(\rho_{\varepsilon}(t))\right)
\phi_{\varepsilon^{-1}}
\left(x_i-m_{\beta,i}^L
(\rho_{\varepsilon}(t))\right).
\]
For almost every $t\in(0,T)$, by Lemma \ref{W2estimate}, we have
\[
m_{\beta,i}^L(\rho_\varepsilon(t))
\longrightarrow
m_{\beta,i}^L(\rho(t)).
\]
Hence, for each such $t$ and every fixed
$x\in\mathbb R^d$, the sequence $\{x_i-m_{\beta,i}^L(\rho_\varepsilon(t))\}$
is bounded. Since $\varepsilon^{-1}\to\infty$,
the defining properties of the cutoff function
imply that, for sufficiently small $\varepsilon$,
\[
\phi_{\varepsilon^{-1}}
\left(x_i-m_{\beta,i}^L
(\rho_\varepsilon(t))\right)=1.
\]
Consequently,
\begin{equation}
h_i^\varepsilon(t,x)
\longrightarrow
x_i-m_{\beta,i}^L(\rho(t))
=:h_i(t,x)
\quad\text{for a.e. }(t,x)
\in(0,T)\times\mathbb R^d.
\label{eq:exist-h-convergence}
\end{equation}

Moreover, let $\Omega\subset\R^d$ be a bounded
domain and choose $R_\Omega>0$ such that
$\Omega\subset B_{R_\Omega}(0)$.
By Proposition \ref{consensus point prop} with $p=2$ and the
uniform second moment estimates \eqref{M2bound}, we have
\[
|m_\beta^L(\rho_\varepsilon(t))|^2
\le C'\int_{\R^d}|x|^2
\rho_\varepsilon(t,x)\,dx
\le C'M_2.
\]
Similarly, by \eqref{M2bound-limit}, $|m_\beta^L(\rho(t))|^2
\le C'M_2.$
Using $|\phi_R|\le1$, we obtain
\begin{equation}
|h_i^\varepsilon(t,x)|+|h_i(t,x)|
\le 2R_\Omega+2\sqrt{C'M_2}
=:C_{\Omega,T},
\label{eq:exist-h-bound}
\end{equation}
for almost every $t\in(0,T)$ and all
$x\in\Omega$, uniformly with respect
to $\varepsilon$.

We now establish the convergence of
the three terms separately.

\medskip
\noindent\textbf{The time derivative term.}

By the convergence of  $\rho_\varepsilon$ to $\rho$ in $L^1((0,T)\times\Omega)$ in \eqref{eq:globalL1}, we have
\begin{align*}
&\left|
\int_0^T\int_\Omega
(\rho_{\varepsilon}-\rho)
\partial_t\varphi\,dxdt
\right|\le
\|\partial_t\varphi\|_{L^\infty}
\|\rho_{\varepsilon}-\rho\|
_{L^1((0,T)\times\Omega)}
\longrightarrow0\quad\text{as}~\varepsilon\to0.
\end{align*}

\medskip
\noindent\textbf{The drift term.}

For each $i=1,\ldots,d$, we consider
\begin{align}
\left|\int_0^T\int_{\mathbb{R}^d}
(h_i^\varepsilon\rho_\varepsilon-h_i\rho)\,\partial_i\varphi\,dx\,dt\right|
\le{}&
\left|
\int_0^T\int_\Omega
h_i^{\varepsilon}
(\rho_{\varepsilon}-\rho)
\partial_{x_i}\varphi\,dxdt
\right|
+
\left|
\int_0^T\int_\Omega
(h_i^{\varepsilon}-h_i)\rho
\partial_{x_i}\varphi\,dxdt
\right|
\nonumber\\
=:{}&
I_{i,1}
+I_{i,2}
\label{eq:exist-drift-decomposition}
\end{align}

For the first term, by
\eqref{eq:globalL1} and \eqref{eq:exist-h-bound}, we have
\begin{align*}
I_{i,1}
&\le
C_{\Omega,T}
\|\partial_{x_i}\varphi\|_{L^\infty}
\|\rho_{\varepsilon}-\rho\|
_{L^1((0,T)\times\Omega)}\longrightarrow0\quad\text{as}~ \varepsilon\to0.
\end{align*}

For the second term, by
\eqref{eq:exist-h-convergence}, we have
\[
(h_i^{\varepsilon}-h_i)\rho
\partial_{x_i}\varphi
\longrightarrow0
\quad\text{for a.e. }
(t,x)\in(0,T)\times\Omega.
\]
Moreover, by \eqref{eq:exist-h-bound}, $\left|
(h_i^{\varepsilon}-h_i)\rho
\partial_{x_i}\varphi
\right|
\le
C_{\Omega,T}
|\rho|\|\partial_{x_i}\varphi\|_{L^\infty}
.$
Since
$\rho\in L^1((0,T)\times\Omega)$,
the dominated convergence theorem yields
\[
I_{i,2}
\longrightarrow0
\quad\text{as }\varepsilon\to0.
\]

Combining the above estimates,
we conclude that
\begin{equation}
\begin{aligned}
&\int_0^T\int_\Omega
h_i^{\varepsilon}\rho_{\varepsilon}
\partial_{x_i}\varphi\,dxdt\longrightarrow
\int_0^T\int_\Omega
h_i\rho\partial_{x_i}\varphi\,dxdt.
\end{aligned}
\label{eq:exist-drift-limit}
\end{equation}

\medskip
\noindent\textbf{The diffusion term.}

For each $i=1,\ldots,d$, we consider
\begin{align}
\left|
\int_0^T\int_\Omega
\left[
(h_i^{\varepsilon})^2
\rho_{\varepsilon}
-h_i^2\rho
\right]
\partial_{x_ix_i}\varphi\,dxdt
\right|
\le{}&
\left|
\int_0^T\int_\Omega
(h_i^{\varepsilon})^2
(\rho_{\varepsilon}-\rho)
\partial_{x_ix_i}\varphi\,dxdt
\right|\nonumber\\
+&
\left|
\int_0^T\int_\Omega
\left[
(h_i^{\varepsilon})^2-h_i^2
\right]\rho
\partial_{x_ix_i}\varphi\,dxdt
\right|
=:{}
D_{i,1}
+D_{i,2}.
\label{eq:exist-diffusion-decomposition}
\end{align}

For the first term, by \eqref{eq:globalL1} and
\eqref{eq:exist-h-bound}, we obtain
\begin{align*}
D_{i,1}
&\le
C_{\Omega,T}^2
\|\partial_{x_ix_i}\varphi\|_{L^\infty}
\|\rho_{\varepsilon}-\rho\|
_{L^1((0,T)\times\Omega)}\longrightarrow0
\quad\text{as }~\varepsilon\to0.\end{align*}

For the second term, by
\eqref{eq:exist-h-convergence},
we have
\[
(h_i^{\varepsilon})^2-h_i^2
\longrightarrow0
\quad\text{for a.e. }
(t,x)\in(0,T)\times\Omega.
\]
Furthermore, \eqref{eq:exist-h-bound}
implies $\left|
(h_i^{\varepsilon})^2-h_i^2
\right|
\le C_{\Omega,T}^2.$
Hence,
\begin{align*}
&\left|
\left[
(h_i^{\varepsilon})^2-h_i^2
\right]\rho
\partial_{x_ix_i}\varphi
\right|\le
C_{\Omega,T}^2
\|\partial_{x_ix_i}\varphi\|_{L^\infty}
|\rho|.
\end{align*}
Since the right-hand side is integrable
over $(0,T)\times\Omega$, the dominated
convergence theorem gives
\[
D_{i,2}
\longrightarrow0
\quad\text{as }\varepsilon\to0.
\]
Consequently,
summing over $i=1,\ldots,d$ yields
\begin{align}
&\sum_{i=1}^d
\int_0^T\int_\Omega
(h_i^{\varepsilon})^2
\rho_{\varepsilon}
\partial_{x_ix_i}\varphi\,dxdt
\longrightarrow
\sum_{i=1}^d
\int_0^T\int_\Omega
h_i^2\rho
\partial_{x_ix_i}\varphi\,dxdt.
\label{eq:exist-anisotropic-limit}
\end{align}

The artificial ellipticity term vanishes. Indeed, by conservation
of mass,
\begin{align}
\left|
\varepsilon^2
\int_0^T\int_{\mathbb{R}^d}
\rho_\varepsilon\,\partial_{ii}\varphi\,dx\,dt
\right|
&\leq
\varepsilon^2
\|\partial_{ii}\varphi\|_{L^\infty}
\int_0^T\int_{\mathbb{R}^d}\rho_\varepsilon\,dx\,dt
\longrightarrow0.
\label{eq:epsilon-limit}
\end{align}

\medskip
\noindent\textbf{The initial term.}
Since $\rho_{0,\varepsilon}=J_\varepsilon*\rho_0$
and $\rho_0\in L^1(\mathbb R^d)$, we have
$\rho_{0,\varepsilon}\longrightarrow\rho_0
\quad\text{in }L^1(\mathbb R^d).$ Thus,
\[
\int_{\mathbb R^d}\rho_{0,\varepsilon}(x)\varphi(0,x)\,dx
\longrightarrow
\int_{\mathbb R^d}\rho_0(x)\varphi(0,x)\,dx.
\]

Passing to the limit in \eqref{eq:weak-regularized}, we obtain
\begin{align}
-\int_{\mathbb{R}^d}
\rho_0(x)\varphi(0,x)\,dx={}&
\int_0^T\int_{\mathbb{R}^d}
\rho\,\partial_t\varphi\,dx\,dt
-\lambda\sum_{i=1}^d
\int_0^T\int_{\mathbb{R}^d}
\bigl(x_i-m_{\beta,i}^L(\rho(t))\bigr)
\rho\,\partial_i\varphi\,dx\,dt
\nonumber\\
&+\frac{\sigma^2}{2}\sum_{i=1}^d
\int_0^T\int_{\mathbb{R}^d}
\bigl(x_i-m_{\beta,i}^L(\rho(t))\bigr)^2
\rho\,\partial_{ii}\varphi\,dx\,dt.
\label{integral equality}
\end{align}

The following lemma ensures that the limit density admits a unique
narrowly continuous representative attaining the prescribed initial datum.
\begin{lemma}[Narrowly continuous representative]
\label{lem:narrow-representative}
Let $T>0$, let $L$ satisfy Assumption~\ref{ass}, and let $\rho_0$
be a probability density with finite fourth moment.
Suppose that $\rho:(0,T)\times\mathbb R^d\to[0,\infty)$
is measurable and satisfies
\[
\int_{\mathbb R^d}\rho(t,x)\,dx=1
\quad\text{for a.e. }t\in(0,T),
\qquad
\operatorname*{ess\,sup}_{t\in(0,T)}
\int_{\mathbb R^d}|x|^4\rho(t,x)\,dx\leq M_4.
\]
Assume also that $\rho$ satisfies the integral identity in
Definition~\ref{Definition of weak solutions}
with initial datum $\rho_0$. Then there exists a unique narrowly continuous curve
$\{\rho_t\}_{t\in[0,T]}\subset\mathcal P_4(\mathbb R^d)$
such that
\[
\rho_t(dx)=\rho(t,x)\,dx
\quad\text{for a.e. }t\in(0,T).
\]
Moreover,
\[
\rho_{t=0}(dx)=\rho_0(x)\,dx,
\qquad
\sup_{t\in[0,T]}
\int_{\mathbb R^d}|x|^4\,\rho_t(dx)\leq M_4.
\]
\end{lemma}

\begin{proof}
For $\varphi\in C_c^\infty(\mathbb R^d)$, define
$A_\varphi(t):=\int_{\mathbb R^d}\varphi(x)\rho(t,x)\,dx$
for almost every $t\in(0,T)$. For any $\eta\in C_c^\infty((0,T))$,
we test \eqref{integral equality} with $\psi(t,x)=\eta(t)\varphi(x)$.
Since $\eta$ vanishes near $t=0$, the initial term vanishes, and
the definition of the distributional derivative gives
%For $\varphi\in C_c^\infty(\mathbb R^d)$, define
%\[
%A_\varphi(t):=
%$\int_{\mathbb R^d}\varphi(x)\rho(t,x)\,dx
%\]
%for almost every $t\in(0,T)$.
%To identify the distributional time derivative of $A_\varphi$,
%let $\eta\in C_c^\infty((0,T))$ be arbitrary and choose 
%\psi(t,x):=\eta(t)\varphi(x)$
%as the test function in the integral identity. 
%Since $\eta$ vanishes near $t=0$, the initial term is zero.
%Substituting $\psi$ into the integral identity \eqref{integral equality} and using the
%definition of $A_\varphi$, we obtain
%\[
%\begin{aligned}
%-\int_0^T A_\varphi(t)\eta'(t)\,dt
%={}&-\lambda\sum_{i=1}^d
%\int_0^T\eta(t)
%\left[
%\int_{\mathbb R^d}
%\bigl(x_i-m_{\beta,i}^L(\rho(t))\bigr)
%\rho(t,x)\partial_i\varphi(x)\,dx
%\right]dt
%\\
%&+\frac{\sigma^2}{2}\sum_{i=1}^d
%\int_0^T\eta(t)
%\left[
%\int_{\mathbb R^d}
%\bigl(x_i-m_{\beta,i}^L(\rho(t))\bigr)^2
%\rho(t,x)\partial_{ii}\varphi(x)\,dx
%\right]dt.
%\end{aligned}
%\]
%As this identity holds for every
%$\eta\in C_c^\infty((0,T))$, the definition of the
%distributional derivative yields
\[
\begin{aligned}
A_\varphi'(t)
={}&-\lambda\sum_{i=1}^d
\int_{\mathbb R^d}
\bigl(x_i-m_{\beta,i}^L(\rho(t))\bigr)
\rho(t,x)\partial_i\varphi(x)\,dx
+\frac{\sigma^2}{2}\sum_{i=1}^d
\int_{\mathbb R^d}
\bigl(x_i-m_{\beta,i}^L(\rho(t))\bigr)^2
\rho(t,x)\partial_{ii}\varphi(x)\,dx
\end{aligned}
\]
in the distributional sense on $(0,T)$. By Proposition~\ref{consensus point prop} and the moment bound, $m_{\beta,i}^L(\rho(t))$ is
essentially bounded on $(0,T)$.
Since $\varphi$ is compactly supported and $\rho(t)$ has
unit mass, the right-hand side belongs to $L^\infty(0,T)$.
Thus $A_\varphi$ admits a Lipschitz continuous representative
on $[0,T]$, still denoted by $A_\varphi$.
The initial term in the integral identity yields $A_\varphi(0)=\int_{\mathbb R^d}\varphi(x)\rho_0(x)\,dx.$

Choose $\{\varphi_j\}_{j\ge1}\subset C_c^\infty(\mathbb R^d)$
dense in $C_0(\mathbb R^d)$ under the uniform norm.
Define
\[
E:=
\left\{
s\in(0,T):
\begin{array}{l}
\displaystyle
\int_{\mathbb R^d}\rho(s,x)\,dx=1,\int_{\mathbb R^d}|x|^4\rho(s,x)\,dx\le M_4,\\[6pt]
\displaystyle
\int_{\mathbb R^d}\varphi_j(x)\rho(s,x)\,dx
=A_{\varphi_j}(s)
\quad\text{for every }j\ge1
\end{array}
\right\}.
\]
Since these conditions hold almost everywhere and the family
$\{\varphi_j\}_{j\ge1}$ is countable, $E$ has full measure
in $(0,T)$.

For $s\in E$, set $\mu_s(dx):=\rho(s,x)\,dx$. Fix $t\in[0,T]$ and choose $s_n\in E$ with $s_n\to t$.
The estimate
\[
\mu_{s_n}(\{|x|>R\})
\le \frac1{R^4}\int|x|^4\,d\mu_{s_n}
\le\frac{M_4}{R^4}
\]
implies tightness.
Hence, by Prokhorov's theorem, there exist a subsequence
$\{\mu_{s_{n_k}}\}$ and a probability measure $\mu$ such that
$\mu_{s_{n_k}}$ converges narrowly to $\mu$.
For every $j\ge1$, narrow convergence and the continuity of
$A_{\varphi_j}$ yield
\[
\int_{\mathbb R^d}\varphi_j\,d\mu
=
\lim_{k\to\infty}
\int_{\mathbb R^d}\varphi_j\,d\mu_{s_{n_k}}
=
\lim_{k\to\infty}A_{\varphi_j}(s_{n_k})
=
A_{\varphi_j}(t).
\]
Since $\{\varphi_j\}$ is dense in $C_0(\mathbb R^d)$,
these identities uniquely determine $\mu$.

We denote this measure by $\rho_t$.
It is independent of the chosen sequence and subsequence,
and
\[
\int_{\mathbb R^d}\varphi_j\,d\rho_t
=A_{\varphi_j}(t)
\quad\text{for all }t\in[0,T],\ j\geq1.
\]
For $t\in E$, the same identities hold for $\mu_t$;
therefore $\rho_t=\mu_t$.
Moreover, lower semicontinuity gives
\begin{align}\label{rhot}
\int_{\mathbb R^d}|x|^4\,d\rho_t
\leq
\liminf_{k\to\infty}
\int_{\mathbb R^d}|x|^4\,d\mu_{s_{n_k}}
\leq M_4.   
\end{align}

To prove narrow continuity, let $t_n\to t$ in $[0,T]$. For every $R>0$, the moment estimate \eqref{rhot} gives 
\[
\rho_{t_n}(\{|x|>R\})
\le
\frac{1}{R^4}
\int_{\mathbb R^d}|x|^4\,d\rho_{t_n}
\le
\frac{M_4}{R^4},
\qquad n\ge1.
\]
Since $\overline B_R$ is compact and $M_4/R^4\to0$ as
$R\to\infty$, the sequence $\{\rho_{t_n}\}$ is tight. By Prokhorov's theorem, every subsequence has a further
subsequence converging narrowly to a probability measure.
Every such limit $\nu$ satisfies
\[
\int_{\mathbb R^d}\varphi_j\,d\nu
=
\lim_{k\to\infty}A_{\varphi_j}(t_{n_k})
=
A_{\varphi_j}(t)
=
\int_{\mathbb R^d}\varphi_j\,d\rho_t
\]
for all $j\geq1$.
Thus $\nu=\rho_t$, and tightness implies that the entire
sequence converges narrowly to $\rho_t$.
Consequently, $t\mapsto\rho_t$ is narrowly continuous.
The identities at $t=0$ imply
$\rho_{t=0}(dx)=\rho_0(x)\,dx$.

Finally, any other narrowly continuous representative has
the same continuous test-function pairings almost everywhere,
and hence everywhere.
The determining property of $\{\varphi_j\}$ then gives
uniqueness.
\end{proof}
By nonnegativity of the
approximate solutions, we have $\rho(t,x)\geq0$
for almost every $(t,x)\in(0,T)\times\mathbb R^d$.
The fourth-moment
bound follows from Fatou's lemma, giving
\[
\int_{\mathbb R^d}|x|^4\rho(t,x)\,dx
\leq
\liminf_{\varepsilon\to0}
\int_{\mathbb R^d}|x|^4\rho_{\varepsilon}(t,x)\,dx
\leq M_4
\]
for almost every $t\in(0,T)$. Applying Lemma~\ref{lem:narrow-representative} to the
limiting integral identity, we obtain a narrowly continuous
representative $\rho_t$ with initial value $\rho_0(x)\,dx$.
Together with the previously established regularity,
this proves that $\rho$ is a weak solution in the sense of
Definition~\ref{Definition of weak solutions}.

Having established the existence of weak solutions,
we now turn to their uniqueness.

\begin{theorem}[Uniqueness]
Assume that $
L$ satisfies Assumption~\ref{ass}, and let the non-negative
initial data $\rho_0\in L^1(\mathbb{R}^d)\cap L^\infty(\mathbb{R}^d)
\cap\mathcal{P}_4(\mathbb{R}^d).$
Then the weak solution $\rho$ to \eqref{component CBO} is unique.
\end{theorem}

\begin{proof}
Let $\rho$ be any weak solution to the nonlinear PDE
\eqref{component CBO}.
We consider the following linear SDE
\begin{equation}
dX_t
=
-\lambda\bigl(X_t-m_\beta^L(\rho_t)\bigr)\,dt
+
\sigma D\bigl(X_t-m_\beta^L(\rho_t)\bigr)\,dB_t,
\label{linear-sde}
\end{equation}
where $D(z)=\operatorname{diag}(z_1,\ldots,z_d).$
By the definition of the solution class,
the curve $t\mapsto\rho_t$ is narrowly
continuous on $[0,T]$.
Accordingly, the consensus point
$m_\beta^L(\rho_t)$ is well defined
for every $t\in[0,T]$.

By Proposition~\ref{consensus point prop}
and the uniform second moment bound,
we have $|m_\beta^L(\rho_t)|^2
\le CM_2.$
Hence, the drift and diffusion
coefficients satisfy
\[
|b(t,x)|+\|\Sigma(t,x)\|_F
\le C_T(1+|x|),
\]
where $b(t,x)=-\lambda(x-m_\beta^L(\rho_t))$ and $\Sigma(t,x)=\sigma D(x-m_\beta^L(\rho_t))$.
Moreover,
\[
|b(t,x)-b(t,y)|
+\|\Sigma(t,x)-\Sigma(t,y)\|_F
\le(\lambda+\sigma)|x-y|.
\]

For the prescribed curve $m_\beta^L(\rho_t)$,
the coefficients in \eqref{linear-sde} are
globally Lipschitz with respect to $x$ and
satisfy a linear growth condition.
Hence, by standard SDE theory,
\eqref{linear-sde} admits a unique strong
solution $X$ with initial law $\rho_0$. 
Define $\rho_t^X:=\operatorname{Law}(X_t).$
By It\^o's formula, $\rho^X$ is a weak solution
to the linear Fokker--Planck equation
\begin{equation}
\partial_t\rho^X
=
\lambda\nabla\cdot
\bigl((x-m_\beta^L(\rho_t))\rho^X\bigr)
+
\frac{\sigma^2}{2}
\sum_{i=1}^d
\partial_{ii}
\left(
(x_i-m_{\beta,i}^L(\rho_t))^2\rho^X
\right).
\label{linear-pde}
\end{equation}

We next show that the solution to
\eqref{linear-pde} is unique.
Let $\rho_1^X$ and $\rho_2^X$ be two solutions
to \eqref{linear-pde} with the same initial
data, and set $\bar\rho^X=\rho_1^X-\rho_2^X.$
Then $\bar\rho^X$ satisfies
\[
\partial_t\bar\rho^X
=
\lambda\nabla\cdot
\bigl((x-m_\beta^L(\rho_t))\bar\rho^X\bigr)
+
\frac{\sigma^2}{2}
\sum_{i=1}^d
\partial_{ii}
\left(
(x_i-m_{\beta,i}^L(\rho_t))^2\bar\rho^X
\right),
\qquad
\bar\rho^X(0,x)=0.
\]

For sufficiently regular solutions,
multiplying the above equation by
$\bar\rho^X$ and integrating over
$\mathbb{R}^d$, we obtain
\begin{align*}
\frac12\frac{d}{dt}
\|\bar\rho^X\|_{L^2(\R^d)}^2
&=
-\lambda\int_{\mathbb{R}^d}
\nabla\bar\rho^X\cdot
(x-m_\beta^L(\rho_t))
\bar\rho^X\,dx
-\frac{\sigma^2}{2}
\sum_{i=1}^d
\int_{\mathbb{R}^d}
\partial_i\bar\rho^X\,
\partial_i\left(
(x_i-m_{\beta,i}^L(\rho_t))^2
\bar\rho^X
\right)\,dx\\
&=
\frac{d}{2}(\lambda+\sigma^2)
\|\bar\rho^X\|_{L^2(\R^d)}^2
-\frac{\sigma^2}{2}
\sum_{i=1}^d
\int_{\mathbb{R}^d}
(x_i-m_{\beta,i}^L(\rho_t))^2
|\partial_i\bar\rho^X|^2\,dx\\
&\le
\frac{d}{2}(\lambda+\sigma^2)
\|\bar\rho^X\|_{L^2(\R^d)}^2.
\end{align*}
Therefore, \text{Gr\"onwall}'s inequality and
$\bar\rho^X(0)=0$ imply $\bar\rho^X(t)=0$ in $L^2(\mathbb{R}^d)$ for almost every $t\in(0,T)$.

For general weak solutions, let $\mu$ be any narrowly
continuous probability-valued solution of the linear
Fokker--Planck equation above, with the prescribed initial
law and uniformly bounded second moments.
The linear growth of the coefficients implies
\[
\begin{aligned}
&\int_0^T\int_{\mathbb R^d}
\bigl(|b(t,x)|+\|\Sigma(t,x)\|_{\mathrm F}^2\bigr)
\,d\mu_t(x)\,dt
\le
C_T\int_0^T
\left(1+\int_{\mathbb R^d}|x|^2\,d\mu_t(x)\right)dt
<\infty.
\end{aligned}
\]
By the superposition principle
\cite[Theorem~2.5]{Trevisan2016},
$\mu$ is the curve of time marginals of a solution
to the associated martingale problem.
Since the frozen-coefficient SDE is well posed in law,
its associated martingale problem is also well posed.
Thus the linear Fokker--Planck equation has at most one
solution in this class.
Both $\rho$ and $\rho^X$ belong to this class.

On the other hand, by construction, $\rho$
itself is also a weak solution to
\eqref{linear-pde}.
Hence, by the uniqueness of the linear equation,
\[
\rho_t=\rho_t^X=\operatorname{Law}(X_t).
\]
Therefore, any weak solution to the nonlinear
PDE \eqref{component CBO} is the distribution
of a solution to the nonlinear McKean--Vlasov SDE
\begin{equation}
dX_t
=
-\lambda\bigl(X_t-m_\beta^L(\rho_t)\bigr)\,dt
+
\sigma D\bigl(X_t-m_\beta^L(\rho_t)\bigr)\,dB_t,
\qquad
\rho_t=\operatorname{Law}(X_t).
\label{nonlinear-sde}
\end{equation}

Let $\rho^X$ and $\rho^Y$ be two weak solutions of
the nonlinear PDE with the same initial datum.
On a common probability space, choose an initial
random variable $X_0$ with law $\rho_0$, independent
of a $d$-dimensional Brownian motion $B$.
Let $X$ and $Y$ solve the frozen-coefficient SDEs
corresponding to the prescribed curves
$m_\beta^L(\rho_t^X)$ and $m_\beta^L(\rho_t^Y)$,
respectively, with $Y_0=X_0$ and driven by the same $B$.
By the preceding identification,
\[
\operatorname{Law}(X_t)=\rho_t^X,
\qquad
\operatorname{Law}(Y_t)=\rho_t^Y,
\qquad t\in[0,T].
\]
Set $Z_t:=X_t-Y_t,$ and $
a_t:=m_\beta^L(\rho_t^X)-m_\beta^L(\rho_t^Y).$
Then
\[
dZ_t
=
-\lambda(Z_t-a_t)\,dt
+
\sigma D(Z_t-a_t)\,dB_t.
\]

Applying It\^o's formula to $|Z_t|^2$
and taking expectations, we obtain
\begin{align*}
\mathbb{E}|Z_t|^2
&=
-2\lambda\int_0^t
\mathbb{E}[Z_s\cdot(Z_s-a_s)]\,ds
+\sigma^2\int_0^t
\mathbb{E}|Z_s-a_s|^2\,ds,
\end{align*}
where we have used $Z_0=0$.

By Young's inequality,
\[
\mathbb{E}|Z_t|^2
\le
C\int_0^t
\left(
\mathbb{E}|Z_s|^2+|a_s|^2
\right)\,ds.
\]

By Proposition \ref{consensus point prop} and the definition of the $2$-Wasserstein distance (see \cite{villani2009optimal}), there exists $C_T>0$
such that
\begin{align*}
|a_s|^2
=|m_\beta^L(\rho_s^X)-m_\beta^L(\rho_s^Y)|^2&\le C_TW_2^2(\rho_s^X,\rho_s^Y)\le C_T\mathbb E|X_s-Y_s|^2.
\end{align*}
Here $\operatorname{Law}(X_s,Y_s)$ is a coupling
of $\rho_s^X$ and $\rho_s^Y$.
Consequently,
\[
\mathbb{E}|Z_t|^2
\le
C_T\int_0^t\mathbb{E}|Z_s|^2\,ds.
\]
\text{Gr\"onwall}'s inequality and continuity in time imply
\[
\mathbb{E}|Z_t|^2=0,
\qquad
\rho_t^X=\rho_t^Y,
\qquad
\text{for every }t\in[0,T].
\]

Hence, the nonlinear McKean--Vlasov SDE
\eqref{nonlinear-sde} is unique in law.
Since every weak solution to
\eqref{component CBO} corresponds to the
law of a solution to \eqref{nonlinear-sde},
we conclude that the nonlinear PDE
\eqref{component CBO} admits at most one
weak solution with the prescribed initial
datum.
\end{proof}

Next, we proceed to demonstrate the improved $H^2$-regularity property.
\begin{lemma}\textup{($H^2$-regularity)}
Let $L$ satisfy Assumption~\ref{ass}, and let $\rho_0$
be a nonnegative probability density such that $\rho_0\in
L^1(\mathbb R^d)\cap L^\infty(\mathbb R^d)
\cap H^2(\mathbb R^d)\cap\mathcal P_4(\mathbb R^d).$
Then the unique weak solution $\rho$ satisfies
\[
\rho\in L^\infty(0,T;H^2(\mathbb R^d)),
\qquad
\partial_t\rho\in
L^\infty(0,T;L^2_{\mathrm{loc}}(\mathbb R^d)).
\]
\end{lemma}

\begin{proof}
For convenience, write $a_i(t,x):=x_i-\bigl(m_\beta^L(\rho_t)\bigr)_i.$
Since $m_\beta^L(\rho_t)$ is independent of the spatial
variable $x$, we have $\partial_{x_j}a_i=\delta_{ij}$.
Equation \eqref{component CBO} can therefore be written as
\begin{align*}
\partial_t\rho
&=\lambda\sum_{i=1}^d\partial_{x_i}(a_i\rho)
  +\frac{\sigma^2}{2}
   \sum_{i=1}^d\partial_{x_ix_i}(a_i^2\rho)\\
&=d(\lambda+\sigma^2)\rho
  +(\lambda+2\sigma^2)a\cdot D\rho
  +\frac{\sigma^2}{2}
   \sum_{i=1}^d a_i^2\rho_{x_ix_i}.
\end{align*}

First, applying the first-order spatial derivative $D$ to the
equation, multiplying by $D\rho$, and integrating over
$\mathbb{R}^d$, we obtain
\begin{align*}
\frac12\frac{d}{dt}
\int_{\mathbb{R}^d}|D\rho|^2\,dx
&=d(\lambda+\sigma^2)
  \int_{\mathbb{R}^d}|D\rho|^2\,dx
  +J_1^1+J_1^2,
\end{align*}
where
\begin{align*}
J_1^1
:=(\lambda+2\sigma^2)
  \int_{\mathbb{R}^d}
  D\rho\cdot D(a\cdot D\rho)\,dx,\quad
J_1^2
:=\frac{\sigma^2}{2}
  \sum_{i,j=1}^d
  \int_{\mathbb{R}^d}
  \rho_{x_j}\,
  \partial_{x_j}(a_i^2\rho_{x_ix_i})\,dx.
\end{align*}
Integration by parts gives
\begin{align*}
J_1^1
&=(\lambda+2\sigma^2)
  \int_{\mathbb{R}^d}
  \left(
  |D\rho|^2+\frac12a\cdot D|D\rho|^2
  \right)\,dx\\
&=\left(1-\frac d2\right)(\lambda+2\sigma^2)
  \int_{\mathbb{R}^d}|D\rho|^2\,dx.
\end{align*}
For the diffusion term, we have
\begin{align*}
J_1^2
&=\frac{\sigma^2}{2}
  \sum_{i,j=1}^d
  \int_{\mathbb{R}^d}
  \rho_{x_j}
  \left(
  2\delta_{ij}a_i\rho_{x_ix_i}
  +a_i^2\rho_{x_ix_ix_j}
  \right)\,dx\\
&=\frac{\sigma^2(d-1)}{2}
  \int_{\mathbb{R}^d}|D\rho|^2\,dx
  -\frac{\sigma^2}{2}
  \sum_{i,j=1}^d
  \int_{\mathbb{R}^d}
  a_i^2|\rho_{x_ix_j}|^2\,dx.
\end{align*}
Consequently,
\[
\frac12\frac{d}{dt}
\int_{\mathbb{R}^d}|D\rho|^2\,dx
\le C(\lambda,\sigma,d)
\int_{\mathbb{R}^d}|D\rho|^2\,dx.
\]
By \text{Gr\"onwall}'s inequality,
\begin{equation}\label{anisotropic-H1}
\|D\rho\|_{L^\infty(0,T;L^2(\mathbb{R}^d))}
\le C.
\end{equation}

Second, applying $D^2$ to the equation, multiplying by $D^2\rho$,
and integrating over $\mathbb{R}^d$, we obtain
\begin{align*}
\frac12\frac{d}{dt}
\int_{\mathbb{R}^d}|D^2\rho|^2\,dx
&=d(\lambda+\sigma^2)
  \int_{\mathbb{R}^d}|D^2\rho|^2\,dx
  +J_2^1+J_2^2,
\end{align*}
where
\begin{align*}
J_2^1
:=(\lambda+2\sigma^2)
  \int_{\mathbb{R}^d}
  D^2\rho\cdot D^2(a\cdot D\rho)\,dx,\quad
J_2^2
:=\frac{\sigma^2}{2}
  \sum_{i,j,k=1}^d
  \int_{\mathbb{R}^d}
  \rho_{x_jx_k}\,
  \partial_{x_jx_k}(a_i^2\rho_{x_ix_i})\,dx.
\end{align*}
Here ``$\cdot$'' denotes the Frobenius inner product.
The drift term satisfies
\begin{align*}
J_2^1
&=(\lambda+2\sigma^2)
  \int_{\mathbb{R}^d}
  \left(
  2|D^2\rho|^2+\frac12a\cdot D|D^2\rho|^2
  \right)\,dx\\
&=\left(2-\frac d2\right)(\lambda+2\sigma^2)
  \int_{\mathbb{R}^d}|D^2\rho|^2\,dx.
\end{align*}

For the diffusion term, set $u_{jk}:=\rho_{x_jx_k}$.
Then
\begin{align*}
\partial_{x_jx_k}
\left(\sum_{i=1}^d a_i^2\rho_{x_ix_i}\right)
&=\sum_{i=1}^d a_i^2\partial_{x_ix_i}u_{jk}
  +2a_j\partial_{x_j}u_{jk}\\
&\quad
  +2a_k\partial_{x_k}u_{jk}
  +2\delta_{jk}u_{jk}.
\end{align*}
Multiplying by $u_{jk}$, summing over $j,k$, and integrating
by parts yields
\begin{align*}
J_2^2
&=\frac{\sigma^2(d-2)}{2}
  \int_{\mathbb{R}^d}|D^2\rho|^2\,dx
  +\sigma^2\sum_{i=1}^d
  \int_{\mathbb{R}^d}|\rho_{x_ix_i}|^2\,dx
  -\frac{\sigma^2}{2}
  \sum_{i,j,k=1}^d
  \int_{\mathbb{R}^d}
  a_i^2|\rho_{x_ix_jx_k}|^2\,dx.
\end{align*}
Since $\sum_{i=1}^d|\rho_{x_ix_i}|^2\le |D^2\rho|^2,$
we conclude that
\[
\frac12\frac{d}{dt}
\int_{\mathbb{R}^d}|D^2\rho|^2\,dx
\le C(\lambda,\sigma,d)
\int_{\mathbb{R}^d}|D^2\rho|^2\,dx.
\]
Again, \text{Gr\"onwall}'s inequality gives
\begin{equation}\label{anisotropic-H2}
\|D^2\rho\|_{L^\infty(0,T;L^2(\mathbb{R}^d))}
\le C.
\end{equation}
Together with the previously established $L^2$-estimate and
\eqref{anisotropic-H1}, this proves
\[
\rho\in L^\infty(0,T;H^2(\mathbb{R}^d)).
\]

Finally, the moment estimate implies that
$m_\beta^L(\rho_t)$ is uniformly bounded on $[0,T]$.
Thus, for every bounded open set
$\Omega\subset\mathbb{R}^d$, the expanded equation gives
\begin{align*}
\|\partial_t\rho\|_{L^2(\Omega)}
&\le C_{\Omega,T}
\left(
\|\rho\|_{L^2(\mathbb{R}^d)}
+\|D\rho\|_{L^2(\mathbb{R}^d)}
+\|D^2\rho\|_{L^2(\mathbb{R}^d)}
\right).
\end{align*}
Taking the essential supremum in time, we obtain
\[
\partial_t\rho
\in L^\infty(0,T;L^2_{\mathrm{loc}}(\mathbb{R}^d)).
\]
This completes the proof.
\end{proof}

%\begin{remark}[Compactness mechanism] 
%The compactness argument can be summarized as 
%\[ \boxed{ \begin{array}{c}
%\text{weighted anisotropic gradient estimate}
%+
%L^p,\ p>2
%\\[1mm]
%\Downarrow
%\\
%\text{uniform spatial translations in }L^2_{\mathrm{loc}}
%\\[2mm]
%+
%\\[-1mm]
%\text{negative Sobolev estimate for }\partial_t\rho_\varepsilon
%\\[1mm]
%\Downarrow
%\\
%\text{uniform temporal translations in }L^2_{\mathrm{loc}}
%\\[2mm]
%\Downarrow
%\\
%\text{Kolmogorov--Riesz}
%\\[1mm]
%\Downarrow
%\\
%\rho_{\varepsilon_k}\to\rho
%\text{ strongly in }L^2_{\mathrm{loc}}.
%\end{array}}
%\]
%Thus the role played by a local Sobolev estimate and an Aubin--Lions argument in the isotropic setting is replaced here by direct space--time translation estimates and the Fr\'echet--Kolmogorov theorem.
%\end{remark}

\bibliographystyle{plain}
\bibliography{componentwise.bib}

\end{document}